\documentclass[1p]{elsarticle}
\pdfoutput=1
\journal{Journal of Computational Physics}
\usepackage{lineno,hyperref}
\modulolinenumbers[5]

\usepackage{latexsym, graphicx, epsfig}
\usepackage{subfig}
\usepackage{url}%,hyperref}
\usepackage{xcolor}
\usepackage{array}
\usepackage{amsfonts,amsmath,amssymb}
\usepackage{verbatim}
\usepackage{mathabx}
\usepackage[ruled]{algorithm2e} % must use algorithm2e and not algorithm % ruled gives horizontal lines
\usepackage{ifthen}
\usepackage{booktabs}
\usepackage{colortbl}

\def\bc{\mathbf{c}}

\def\be{\mathbf{e}}
\def\boldf{\mathbf{f}}

\def\bx{\mathbf{x}}

\def\bphi{\boldsymbol{\phi}}
\def\bPhi{\boldsymbol{\Phi}}

\def\bgamma{\boldsymbol{\gamma}}
\def\blambda{\boldsymbol{\lambda}}

\def\reals{\mathbb{R}}

\newcommand{\norm}[2]{\left\lVert #1 \right\rVert_{#2}}
\newcommand{\inprod}[3]{\left( #1, #2 \right)_{#3}}
\newcommand{\abs}[1]{\left\lvert #1 \right\rvert}
\newcommand{\card}[1]{\lvert#1\rvert}

\DeclareMathOperator*{\argmin}{arg\,min}
\newcommand{\subs}[3]{#1_{#2},\ldots,#1_{#3}}

\def\vand{\boldsymbol{\Phi}}
\def\basis{\phi}
\def\col{\tilde{\basis}}
\def\coef{\boldsymbol{\alpha}}
\def\rhs{\mathbf{f}}

\def\stoptol{\varepsilon}

\def\rv{\xi}           %random variable
\def\brv{{\boldsymbol{\xi}}} %random vector
\def\model{f(\brv)}
\def\sur{\hat{f}(\brv)}
\def\dom{{I_{\brv}}}
\def\bXi{\boldsymbol{\Xi}}
\def\cverr{e_{\mathrm{cv}}}

\def\bW{\mathbf{W}}
\def\btheta{\boldsymbol{\theta}}

\def\bset{\Lambda} % the index defining the PCE truncation

\def\RHiLi{\leavevmode\rlap{\hbox to \hsize{\color{blue!20}\leaders\hrule height .8\baselineskip depth .5ex\hfill}}}
\def\EHiLi{\leavevmode\rlap{\hbox to \hsize{\color{green!20}\leaders\hrule height .8\baselineskip depth .5ex\hfill}}}
\def\IHiLi{\leavevmode\rlap{\hbox to \hsize{\color{yellow!20}\leaders\hrule height .8\baselineskip depth .5ex\hfill}}}
\def\SHiLi{\leavevmode\rlap{\hbox to \hsize{\color{red!20}\leaders\hrule height .8\baselineskip depth .5ex\hfill}}}

\begin{document}

% \maketitle
\begin{frontmatter}
\title{Enhancing $\ell_1$-minimization estimates of polynomial chaos expansions using basis selection}
\author[sandia1]{Jakeman, J.D.\fnref{fn1}\corref{cor1}}\ead{jdjakem@sandia.gov}
\author[sandia1]{Eldred, M.S\fnref{fn1}}
\author[sandia1]{Sargsyan, K.\fnref{fn1}}
\address[sandia1]{Sandia National Laboratories, 
		Albuquerque, NM 87185, 
		United States}
\address[sandia2]{Sandia National Laboratories, 
		Livermore, CA 94550, 
		United States}
\cortext[cor1]{Corresponding author}
\fntext[fn1]{Sandia National Laboratories is a multi-program laboratory 
managed and operated by Sandia Corporation, a wholly owned subsidiary of Lockheed 
Martin Corporation, for the U.S. Department of Energy’s National Nuclear Security 
Administration under contract DE-AC04-94AL85000.}
\begin{abstract}
In this paper we present a basis selection method that can be used with
$\ell_1$-minimization to adaptively determine the large coefficients of 
polynomial chaos expansions (PCE). The adaptive construction produces anisotropic
basis sets that have more terms in important dimensions and limits the number of 
unimportant terms that increase mutual coherence and thus degrade the performance of $\ell_1$-minimization. 
The important features and the accuracy of basis selection are demonstrated with a number of numerical examples.
Specifically, we show that for a given computational budget, basis selection produces a more 
accurate PCE than would be obtained if the basis is fixed a priori.
We also demonstrate that basis selection can be applied with non-uniform random variables and 
can leverage gradient information.
\end{abstract}
\end{frontmatter}
% \linenumbers

\section{Introduction}
Quantifying uncertainty in a computational model is essential to building the confidence of stakeholders in the predictions of that model. 
Sources of uncertainty in model predictions can be broadly grouped into two classes, uncertainty arising from model structure and uncertainty arising from the model parameterization.
The effect of these uncertainties must be traced through the model and the effect on the model output (prediction) needs to be quantified. In this paper we will present a method 
for quantifying parametric uncertainty that utilizes the strengths of Polynomial Chaos Expansions (PCE) and $\ell_1$-minimization.

When the computational cost of a simulation model is large, the most popular and effective means of quantifying parametric uncertainty is to construct
an approximation of the response of the model output to variations in the model input. Once built, this surrogate can be interrogated cheaply, without 
further model evaluations, to obtain statistics of interest such as model output moments and distributions. Within the computational science community, the most widely adopted approximation methods
used for Uncertainty Quantification (UQ) are based on generalized polynomial chaos expansions~\cite{ghanem91,xiu02a}, sparse grid interpolation~\cite{jakeman2013localuq,ma09} and
Gaussian process models~\cite{rasmussen2005}. 

Polynomial chaos expansions represent a response surface as a linear combination of orthonormal multivariate polynomials. The choice of the orthonormal polynomials is related to 
the distribution of the model input variables. Provided sufficient smoothness conditions are met, PCEs exhibit fast convergence -- in some cases even
exponential convergence can be obtained~\cite{Babuska_TZ_SIAMNA_2004,xiu02a}. In this paper we will focus on PCEs as they allow one to leverage the
advantages of $\ell_1$-minimization for computing approximations from limited data.

% A major appeal of PCE is the direct relationship of the polynomial basis and the distribution of the random variables which allows the mean and variance 
% and parameter sensitivities, such as Sobol indices and derivatives, to be computed analytically from the expansion.
% The utility of a PCE is determined by the capacity to accurately estimate the expansion coefficients. 

The stochastic Galerkin~\cite{ghanem91,xiu02a} 
and stochastic collocation~\cite{babuska07,mathelin05,tatang97,xiu05} methods are the two main approaches for approximating the PCE coefficients. 
The former is intrusive and so is only feasible when one has the ability to modify the code used to solve the governing
equations of the model. Stochastic collocation, however,  is a non-intrusive sampling based approach that allows the computational model to be treated as a black box.
In this paper we focus on stochastic collocation which involves running the computational model with a set of realizations of the random parameters
and constructing an approximation of corresponding model output. 

Pseudo-spectral projection~\cite{conrad2013,constantine2012}, sparse grid interpolation~\cite{gana07,jakeman2013localuq,ma09,nobile08b}, 
probabilistic multi-element methods~\cite{foo10} are stochastic collocation methods which have been used effectively in many situations. These methods, however,
all require structured samples and/or the ability to iteratively determine the collocation points. 
% Moreover these methods typically require a set of initial samples, the number of which increases rapidly with dimension. 
% \ks{Not sure I understand the arguments in the last two sentences. What we do is free of these constraints? E.g., Least-squares can also do well with unstructured samples, but the strength of CS is not  because it deals with any sampling, but because it selects basis.}

Recently $\ell_1$-minimization has been shown to be an effective method for approximating PCE coefficients from small number of and possibly arbitrarily 
positioned collocation nodes~\cite{blatman2011,doostan2011,mathelin2012,Sargsyan_SNDRT_IJUQ_2014,YangKarniadakis2013}. These methods are very effective when the number of non-zero terms
in the PCE approximation of the model output is small (i.e. sparse) or the magnitude of the PCE coefficients decay rapidly (i.e. compressible). 

The efficacy of $\ell_1$-minimization when used to estimate PCE coefficients is dependent on the rate of the decay of the PCE coefficients, the characteristics of the stochastic collocation samples 
and the truncation of the PCE. The decay of the coefficients is a property of the model and cannot be adjusted to enhance $\ell_1$-recovery. However, the truncation of the PCE and
the sampling of the model inputs can both be controlled. 

Recently some attention has been given to designing sampling strategies to increase the accuracy of sparse PCE~\cite{Rauhut_W_JAT_2012,Xu_Z_2014_SISC,Yan_GX_IJUQ_2012}. 
Almost no attention, however, has been given to the effect of the PCE truncation when using $\ell_1$-minimization. Typically, when using $\ell_1$-minimization,
a total degree truncation is applied to PCE. However the number of terms in this basis grows factorially with the number of model parameters. 
This fast growth in the number of basis terms significantly affects the ability of $\ell_1$-minimization to accurately approximate PCE coefficients. To reduce the growth of a PCE basis
in high dimensions a hyperbolic cross PCE truncation can be employed~\cite{blatman2011}. However, despite the slower growth of the hyperbolic truncation it can perform poorly when the
`true' PCE has large coefficients associated with interaction basis terms.

The goal of this paper is to present a basis selection algorithm that adaptively determines a set of PCE basis terms that enable accurate approximation of PCE coefficients
using $\ell_1$-minimization.
Specifically, we aim to:
\begin{itemize}
 \item Present an iterative algorithm for selecting a polynomial chaos basis that, for a given computational budget, produces a more accurate PCE than would 
 be obtained if the basis is fixed a priori.
 \item Demonstrate numerically that in high dimensions, for which high-order total-degree PCE bases are infeasible, basis selection allows the accurate identification
 of high-order terms that cannot be captured by a low-order total-degree basis.
 \item Demonstrate numerically that even for lower dimensional problems, for which high-order total-degree PCE bases are feasible, 
 basis selection still produces more accurate results than a priori fixed basis sets.
 \item Show that basis selection can leverage function gradients, that for a given computational budget, will produce more accurate 
 approximations than an approximation based solely on function values.
 \item Illustrate that basis selection can be applied with non-uniform random variables.
\end{itemize}
The remainder of this paper is organized as follows: Section~\ref{sec:pce} provides a brief summary of PCEs; 
Section~\ref{sec:compressed-sensing} discuses how to use $\ell_1$-minimization for building a PCE and the need to move away
from a priori-fixed PCE truncations in higher dimensions; Section~\ref{sec:iterative-basis-selection} proposes a new method for iteratively defining PCE truncations; the properties and effectiveness 
of the proposed method are demonstrated numerically in Section~\ref{sec:results}; and conclusions are presented in Section~\ref{sec:conclusions}.
\section{Polynomial chaos expansions}\label{sec:pce}
% Polynomial Chaos methods represent a function $\model\in L_2(\rho(\brv))$ as an expansion of orthonormal polynomials
% \begin{equation}
% \label{eq:pce-integer-index}
% \model\approx\sur=\sum_{i=1}^N\alpha_{i}\phi_{i}(\brv).
% \end{equation}
% We refer to~\eqref{eq:pce-integer-index} as a polynomial chaos expansion (PCE). The PCE basis functions $\{\phi_{i}(\brv)\}$ are tensor products of orthonormal polynomials 
% which are chosen to be 
% orthonormal to a distribution $\rho(\brv)$ of the random vector $\brv$.
% That is 
% \[
% \inprod{\phi_{i}(\brv)}{\phi_{j}(\brv)}{} = \int_{\dom} \phi_{i}(\brv)\phi_{j}(\brv)\rho(\brv) = \delta_{ij}
% \]
% where $\dom$ is the range of the random variables.
% 
% For any function (model output) with finite variance, i.e. $\model\in L_2(\rho(\brv))$,
% a PCE will converge to that function as the number of terms tends to infinity. The rate of convergence is dependent on the regularity of the response surface. If $\model$ 
% is analytical with respect to the random variables then~\eqref{eq:pce-integer-index} converges exponentially in $L_2(\rho(\brv))$-sense~\cite{bieri09}.

Polynomial Chaos methods represent both the model inputs $\btheta=(\theta_1,\ldots,\theta_{\tilde{d}})$ and model output $f(\btheta)$ as an expansion of orthonormal polynomials of random variables 
$\brv=(\rv_1,\ldots,\rv_d)$. Specifically we represent the random inputs 
 as
\begin{equation}
\label{eq:pce-rv-integer-index}
\theta_n\approx\sum_{i=1}^{N_{\theta_n}}\beta_{i}\phi_{i}(\brv),\quad n=1,\ldots,\tilde{d}
\end{equation}
and the model output as
\begin{equation}
\label{eq:pce-integer-index}
f(\btheta(\brv))\approx\sur=\sum_{i=1}^N\alpha_{i}\phi_{i}(\brv).
\end{equation}
We refer to~\eqref{eq:pce-rv-integer-index} and~\eqref{eq:pce-integer-index} as a polynomial chaos expansion (PCE). 
The PCE basis functions $\{\phi_{i}(\brv)\}$ are tensor products of orthonormal polynomials 
which are chosen to be 
orthonormal with respect to the distribution $\rho(\brv)$ of the random vector $\brv$.
That is 
\[
\inprod{\phi_{i}(\brv)}{\phi_{j}(\brv)}{} = \int_{\dom} \phi_{i}(\brv)\phi_{j}(\brv)\rho(\brv) d\brv= \delta_{ij}
\]
where $\dom$ is the range of the random variables.

The random variable (germ) $\brv$ of the PCE is typically related to the distribution of the input variables. For example, if the one-dimensional input variable $\theta$ is uniform on $[a,b]$ then 
$\rv$ is also chosen to be uniform on [-1,1] and $\phi$ are chosen to be Legendre polynomials such that $\theta=\beta_1 + \beta_2\xi=(b+a)/2+\rv (b-a)/2$. For simplicity and without loss of
generality, we will assume that $\brv$ has the same distribution as $\btheta$ and thus we can use the two variables interchangeably (up to a linear transformation which we will ignore).

% For any function (model output) with finite variance, i.e. $\model\in L_2(\rho(\brv))$, a PCE will converge to that function as the number of terms tends to infinity. 
The rate of convergence is dependent on the regularity of the response surface. If $\model$ 
is analytical with respect to the random variables then~\eqref{eq:pce-integer-index} converges exponentially in $L_2(\rho(\brv))$-sense~\cite{bieri09}.

In practice the PCE~\eqref{eq:pce-integer-index} must be truncated. The most common approach is to set a degree $p$ and retain only the multivariate polynomials of degree at most $p$. 
Rewriting~\eqref{eq:pce-integer-index} using the typical multi-dimensional index notation
\begin{equation}
\label{eq:pce-multi-index}
\model\approx\sur=\sum_{\blambda\in\bset}\alpha_{\blambda}\phi_{\blambda}(\brv)
\end{equation}
the total degree basis of degree $p$ is given by
\begin{equation}\label{eq:hyperbolic-index-set} 
 \bset = \bset^d_{p,q} = \{\phi_{\blambda} : \norm{\blambda}{q} \le p\},\quad \blambda = (\lambda_1,\ldots,\lambda_d)
\end{equation}
with $q=1$. The number of terms in this total degree basis
\[
 \text{card}\; \bset^d_{p,1} \equiv P = { d+p \choose d }
\]
grows factorially with dimension. 
This rapid growth limits the applicability of the total degree basis to moderate dimensions or low degree polynomials in higher dimensions.

\begin{figure}[ht]
\centering
\includegraphics[width=0.475\textwidth]{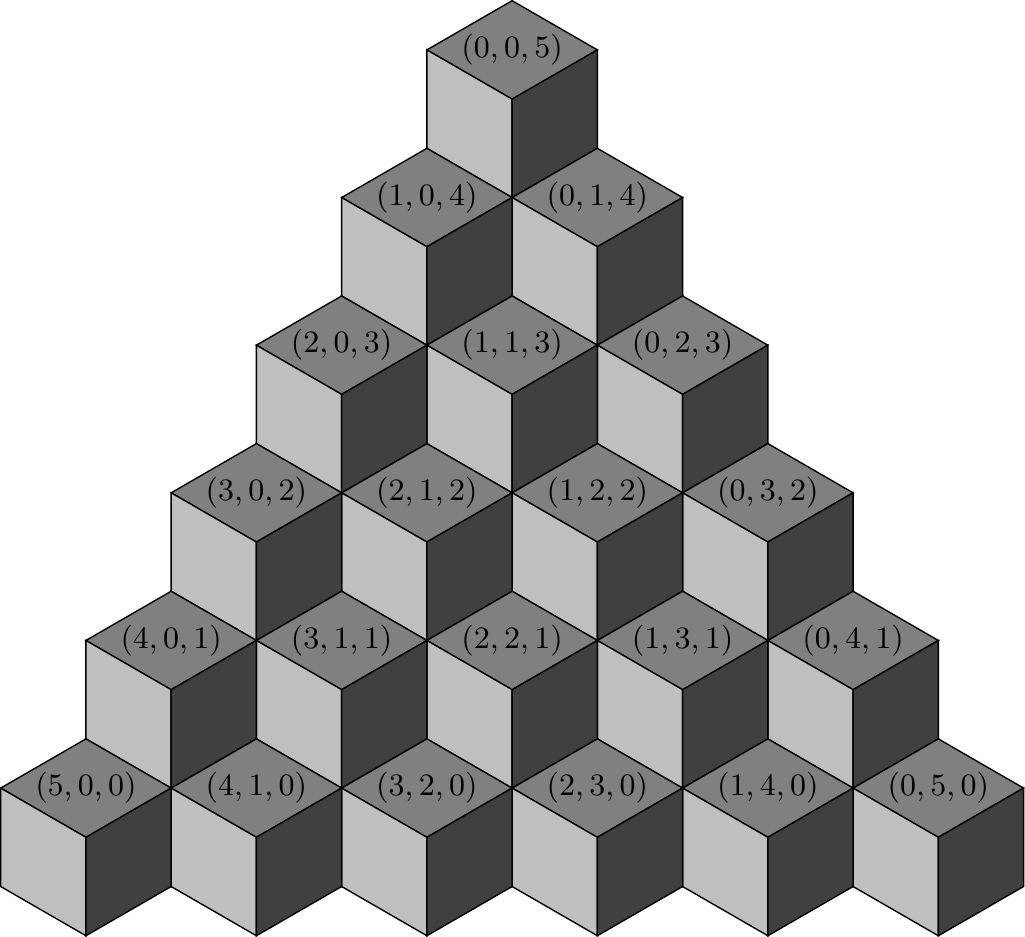}
\includegraphics[width=0.475\textwidth]{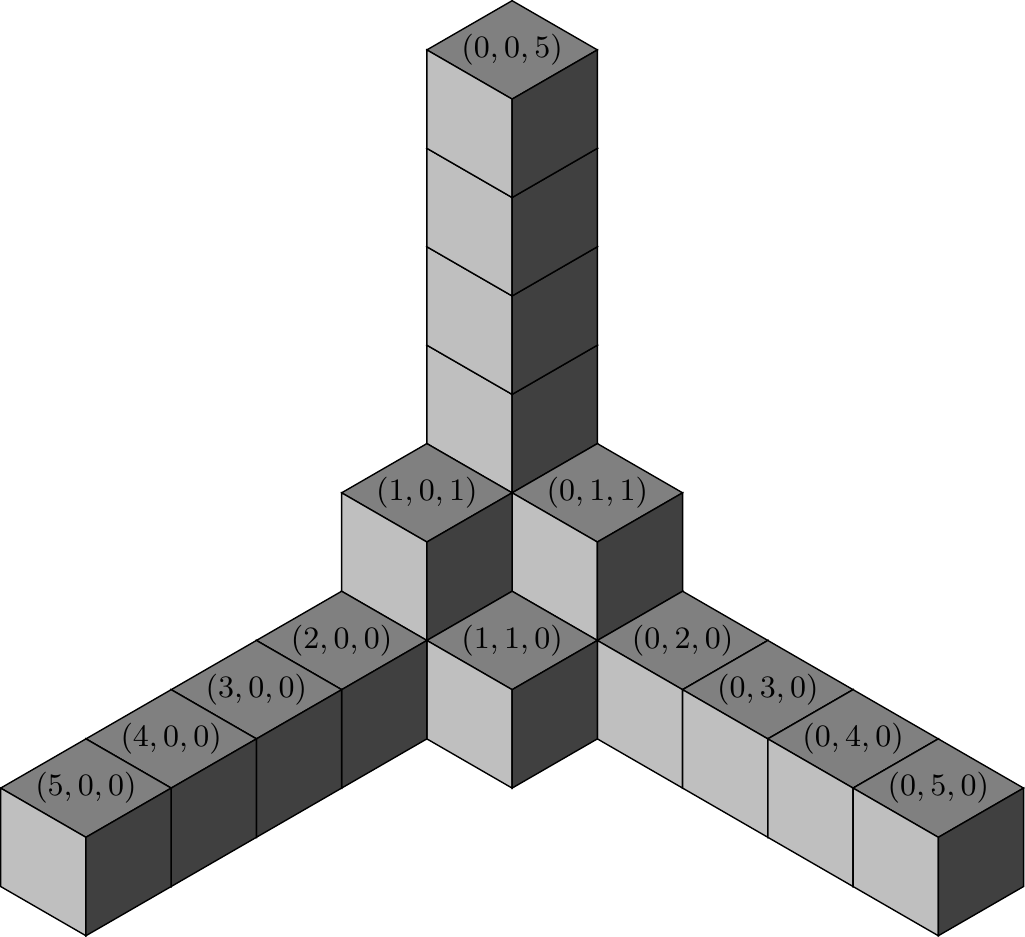}
\caption{(Left) A total degree index set $\bset^3_{6,1}$. (Right) A hyperbolic index set $\bset^3_{6,1/2}$. Each cube represents a $3$-dimensional index $\blambda=(\lambda_1,\lambda_2,\lambda_3)$ 
in $\bset$. The value of each $\blambda$ is given on the top of each cube.}
\label{fig:index-set-comparison}
\end{figure}

The authors of~\cite{blatman2011} propose using hyperbolic index sets, \eqref{eq:hyperbolic-index-set} with $q<1$, to slow the growth of the PCE basis with dimensionality.
The use of hyperbolic indices assumes that the contribution to variance from the interaction between the random variables decays rapidly
as the number of variables involved in the interaction increases. 
Figure~\ref{fig:index-set-comparison} shows a three dimensional total degree and hyperbolic index set. It is clear that for a given degree $p$ the hyperbolic index set has many less terms
than the total degree polynomial basis. However the smaller basis size requires omitting polynomial terms interaction terms, 
that is indices $\blambda$ with at least two $\lambda_n>0, n=1,\ldots,d$. When a function has 
large non-zero PCE coefficients corresponding to these missing multivariate basis terms, the hyperbolic index set may be an inappropriate form of truncation. Ideally the basis truncation
should be adapted to the function being approximated.

\section{$\ell_1$-minimization}\label{sec:compressed-sensing}
The coefficients of a polynomial chaos expansion can be approximated effectively using $\ell_1$-minimization. Specifically,
given a small set of $M$ unstructured realizations  $\bXi=\{\brv_1,\ldots,\brv_M\}$, with 
corresponding model outputs $\boldf=(f(\brv_1),\ldots,f(\brv_M))^T$, we would like to find
a solution that satisfies
\[
 \vand\coef \approx \boldf
\]
where $\coef=(\alpha_{\blambda_1},\ldots,\alpha_{\blambda_N})^T$ denotes the vector of PCE coefficients and $\bPhi$ denotes the Vandermonde matrix with entries
$\bPhi_{ij} = \phi_j(\brv_i),\quad i=1,\ldots,M,\; j=1,\ldots,N$. 

When the model $\model$ is high-dimensional and computationally expensive, and non-adaptive basis truncation rules are employed, the number of model simulations that can be generated is much smaller than the number of unknown PCE coefficients, i.e $M\ll N$. 
Under these conditions, finding the PCE coefficients is ill-posed and we must impose some form of regularization to obtain a unique solution.

$\ell_1$-minimization provides a means of identifying {\it sparse} coefficient vectors from a limited amount of simulation data. 
A polynomial chaos expansion is defined as $s$-sparse when $\norm{\coef}{0} \le s$, i.e. the number of non-zero coefficients does not exceed $s$. In practice, not many simulation models will be truly sparse, 
but PCE are often {\it compressible}, that is the magnitude of the coefficients decay rapidly or alternatively most of the PCE variance is concentrated in a few terms. Compressible vectors are well represented by sparse vectors
and thus the coefficients of compressible PCE can also be recovered accurately using $\ell_1$-minimization.

$\ell_1$-minimization attempts to find the dominant PCE coefficients by solving the following optimization problem
\begin{equation}
\label{eq:bpdn}
\coef = \argmin_{\coef}\; \|\coef\|_1\quad \text{such that}\quad \|\bPhi\coef - \boldf\|_2 \le \stoptol
\end{equation}
This $\ell_1$-minimization problem is often referred to as Basis Pursuit Denoising. The problem obtained by setting $\stoptol=0$, to enforce interpolation, is 
termed Basis Pursuit. There is a close connection between~\eqref{eq:bpdn} and Least Absolute Shrinkage Operator (LASSO)~\cite{tibshirani1996} well known in the statistics literature. 
Indeed these problems are equivalent under certain conditions~\cite{donoho2006}.

\subsection{$\ell_1$-minimization algorithms}
Numerous algorithms~\cite{becker2011,candes2006,needell2010,BergFriedlander2008} exist for solving~\eqref{eq:bpdn} which are all stable and accurate under certain well 
defined conditions. In this paper we will use the greedy algorithm Orthogonal Matching Pursuit (OMP)~\cite{chen2001} to 
estimate PCE coefficients. 
%OMP is explained in detail in Section~\ref{sec:omp}.  
OMP requires stronger theoretical conditions than some of its counterparts~\cite{cai2011} but in practice OMP can still obtain
comparable accuracy to these algorithms. In this paper we use OMP because of its fast execution speed which makes 
OMP more amenable to cross validation which can be used to estimate optimal method parameters such as the tolerance $\stoptol$ of~\eqref{eq:bpdn}. 
We remark, however, that the basis selection procedure presented in this paper can be used in conjunction with most $\ell_1$-minimization algorithms.

\subsubsection{Hyper-parameter estimation via cross validation}
Accurately computing the coefficients of a polynomial chaos expansion requires determining a `good' truncation set $\bset$ and specifying the tolerance $\stoptol$ in the Basis Pursuit DeNoising 
problem~\eqref{eq:bpdn}. Cross validation has been shown to be effective at aiding these choices. Specifically cross validation has been used in the past to estimate
 the polynomial degree $p$ of a hyperbolic expansion~\cite{blatman2011} and 
to estimate the tolerance $\stoptol$ of~\eqref{eq:bpdn}~\cite{blatman2011,boufounos2007,doostan2011,mathelin2012,ward2009}\footnote{The choice of $\stoptol$ 
can significantly affect the accuracy of the PCE obtained using~\eqref{eq:bpdn}.
Decreasing $\stoptol$ can lead to over-fitting, whilst higher values of $\stoptol$ can deteriorate the accuracy of the approximation.}.

In this paper we will use $K=10$ fold cross validation to choose the values of sets of hyper-parameters $\bgamma$. The number and type of hyper-parameters
is dependent on the $\ell_1$-minimization method used in conjunction with cross validation.  As an example consider solving in~\eqref{eq:bpdn} using  
an a priori fixed total degree basis $\bset_{p,1}^d$. The hyper-parameters that can be estimated using cross validation are the degree $p$ and
the tolerance $\stoptol$, that is $\bgamma=(p,\stoptol)$.
%for example $\bgamma=(t,\stoptol)$such as the the degree $p$ of a total-degree PCE basis $\bset_{p,q}^d$ and the tolerance $\stoptol$. 

Let $\zeta:\{1,\ldots,M\}\rightarrow\{1,\ldots,K\}$ be an indexing function that determines the partition of the training data. Furthermore let
$\hat{f}^{-\zeta}$ be the PCE approximation built on the data with the $\zeta$ part removed, then the cross validation error is given by
\begin{equation}\label{eq:cv-error}
\cverr(\bgamma)=\frac{1}{M}\sum_{k=1}^K e_{\zeta(k)},\quad e_{\zeta(k)} = \sum_{j\in\zeta(k)} ( y_j - \hat{f}^{-\zeta(k)}(x_j))^2
\end{equation}
To compute $\cverr$ we divide the data pairs $(\Xi,\boldf)$, based upon the randomly chosen partitions $\zeta(k)$, into $K$ sets (folds) of equal size  $(\Xi_k,\boldf_k)$, $k=1,\ldots,K$. 
A PCE $\hat{f}^{-\zeta(k)}$, is then built on the training data $\Xi_{\mathrm{t}}=\Xi \setminus \Xi_k$ 
with the $k$-th fold removed, using the hyper-parameters $\bgamma$. The remaining data $\Xi_{\mathrm{v}}=\Xi_k$ is then used to estimate the prediction error. 
%We denote the number of samples in the training and validation sets $M_t$ and $M_v$ respectively. By construction $M=M_t+M_v$. 
To estimate the hyper-parameters $\bgamma$ we search over a set of possible values for $\bgamma$
and select $\bgamma=\argmin_{\bgamma} \cverr(\bgamma)$.

\begin{figure}[ht]
\centering
\includegraphics[width=0.95\textwidth]{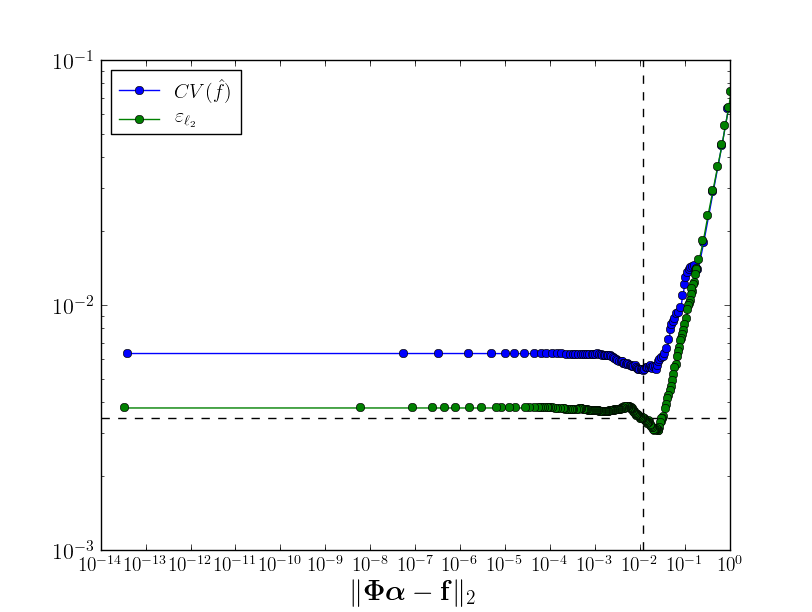}
\caption{The use of cross validation to select the truncation tolerance $\stoptol$ for~\eqref{eq:genz-corner-peak} with coefficients $c^{(1)}$. The vertical line represents the tolerance chosen
by cross validation and the horizontal line is the $\varepsilon_{\ell_2}$ error in the resulting PCE. $M=200$ uniform samples were used. 
Accuracy is measured using the $\ell_2$ norm computed at $100,000$ Latin-hypercube samples (see Section~\ref{sec:results}).}
\label{fig:truncation-tolerance-cross-validation-vs-l2-error}
\end{figure}

Figure~\ref{fig:truncation-tolerance-cross-validation-vs-l2-error}
presents a typical example illustrating the change in the $\ell_2$ error of a PCE with fixed degree, as the tolerance $\stoptol$ is decreased. 
The figure also plots the cross validation error which is a good indicator of the $\ell_2$ error behavior. The vertical line represents the tolerance chosen
by cross validation and the horizontal line is the $\varepsilon_{\ell_2}$ error in the resulting PCE. The result shown is typical. 
There is a bias (underestimation of $\varepsilon_{\ell_2}$) in the cross validation estimate, yet despite this bias cross validation consistently chooses a tolerance that produces a near minimal error.

\subsection{Recoverability of $\ell_1$-minimization}\label{sec:recovery-properties}
The ability of $\ell_1$-minimization to accurately determine the large coefficients of the PCE is determined by the properties of the matrix $\vand$ and the sparsity of PCE representation of 
the model response $\model$. The sparsity is a property of the model and cannot be changed, however the properties of the $\vand$ are influenced by the selection of the 
realizations $\{\brv_i\}_{i=1}^M$ and the truncation $\bset$. 

{\it Mutual coherence} is one measure often used to indicate the ability of $\ell_1$-minimization to find a sparse solution.
The mutual coherence of a matrix $\vand\in\reals^{M\times N}$ with columns $\col_j$ is
\begin{equation}\label{eq:mutual-coherence}
\mu(\vand)= \max_{1<j<k\leq N}\frac{\abs{\col_j^T\col_k}}{\norm{\col_j}{2}\norm{\col_k}{2}}
\end{equation}
and is a measure of the maximum correlation between any two columns in the matrix. $\ell_1$-minimization will obtain a better
estimate of the PCE coefficients if the mutual coherence of $\vand$ is small. Intuitively, if two
columns are closely correlated the mutual coherence will be large and it will be impossible, in general, to distinguish whether the energy in
the signal comes from one or the other. 

The {\it restricted isometry property} (RIP)~\cite{candes2006}, quantified by the restricted isometry constant $\delta$ is another measure of the recoverability of the matrix $\vand$.
For each $s=1,2,\ldots$ the isometry constant $\delta_s$ of a matrix $\vand$ is the smallest number such that
\begin{equation}
\label{eq:rip}
 (1-\delta_s)\norm{\coef_s}{2}^2\le\norm{\vand\coef_s}{2}^2\le(1+\delta_s)\norm{\coef_s}{2}^2
\end{equation}
for all vectors $\coef_s$ with $s$ non-zero entries.
This is equivalent to requiring that the eigenvalues of all Grammian matrices $\vand^\top_{\bset_s}\vand_{\bset_s}$ lie between $[1-\delta_s,1+\delta_s]$, where $\vand_{\bset_s}$ are $M\times s$ submatrices of $\vand$.
The restricted isometry property measures the ability of $\vand$ to preserve the lengths of $s$-sparse vectors. The RIP can be intuitively thought of as a measure of $s$-wise coherence as opposed to 
mutual coherence which is a measure of pair wise coherence.

\subsubsection{Sampling strategies and pre-conditioning}
The sampling strategy used to choose the samples $\Xi$ affects the mutual coherence and RIP of $\vand$ and thus can impact the accuracy of the recovered polynomial chaos expansion. 
To date, the best sampling strategies for $\ell_1$-minimization are random~\cite{Rauhut_W_JAT_2012,Xu_Z_2014_SISC}. The nature of the random samples is dependent on the distribution of the 
random variables $\brv$, the number of random dimensions, and the degree of the PCE. 
%It was shown in~\cite{Yan_GX_IJUQ_2012} that 
%that Legendre based PCE can be recovered with high probability using $M=O(s\log^4N)$ samples chosen randomly from the uniform distribution.

The accuracy of $\ell_1$-minimization solutions of~\eqref{eq:bpdn} can also be improved by the use of pre-conditioning. The pre-conditioned $\ell_1$-minimization problem is given by
\begin{equation}
\label{eq:pre-condition-bpdn}
\coef = \argmin_{\coef}\; \|\coef\|_1\quad \text{such that}\quad \|\bW\bPhi\coef - \bW\boldf\|_2 \le \stoptol
\end{equation}
where $W\in\reals^{M\times M}$ is a diagonal matrix with entries chosen to enhance the recovery properties of $\ell_1$-minimization.
When recovering $s$-sparse one-dimensional Legendre polynomials, randomly sampling $\Xi=\{\xi_m\}_{m=1}^M$ from the Chebyshev measure and 
choosing weights $w_{m,m}=(\pi/2)^{d/2}(1-\xi_m^2)^{1/4}$, can result in significant increases in the accuracy of the coefficients recovered by $\ell_1$-minimization~\cite{Rauhut_W_JAT_2012}.
In the multivariate setting, however, the benefit of pre-conditioning is less clear~\cite{Yan_GX_IJUQ_2012}. 
In this paper all numerical results presented are generated 
without pre-conditioning.

\subsubsection{PCE truncation}
\begin{figure}[ht]
\centering
% \includegraphics[width=0.45\textwidth]{d-6-mutual-coherence-vs-degree-log}
% \includegraphics[width=0.45\textwidth]{d-10-mutual-coherence-vs-degree-log}\\
% \includegraphics[width=0.45\textwidth]{d-6-mutual-coherence-vs-degree-log-random}
% \includegraphics[width=0.45\textwidth]{d-10-mutual-coherence-vs-degree-log-random}
% \caption{The dependence of mutual coherence on the number of terms $P$ in a total-degree multivariate Legendre PCE basis. (Top-left) $d=6$ and (Top-right) $d=10$ using Latin hypercube samples. 
% (Bottom-left) $d=6$ and (Bottom-right) $d=10$ using uniform samples. Each curve is the average of twenty different point designs.}
\makebox[1\linewidth][c]{%
\subfloat[]{\includegraphics[width=0.7\textwidth]{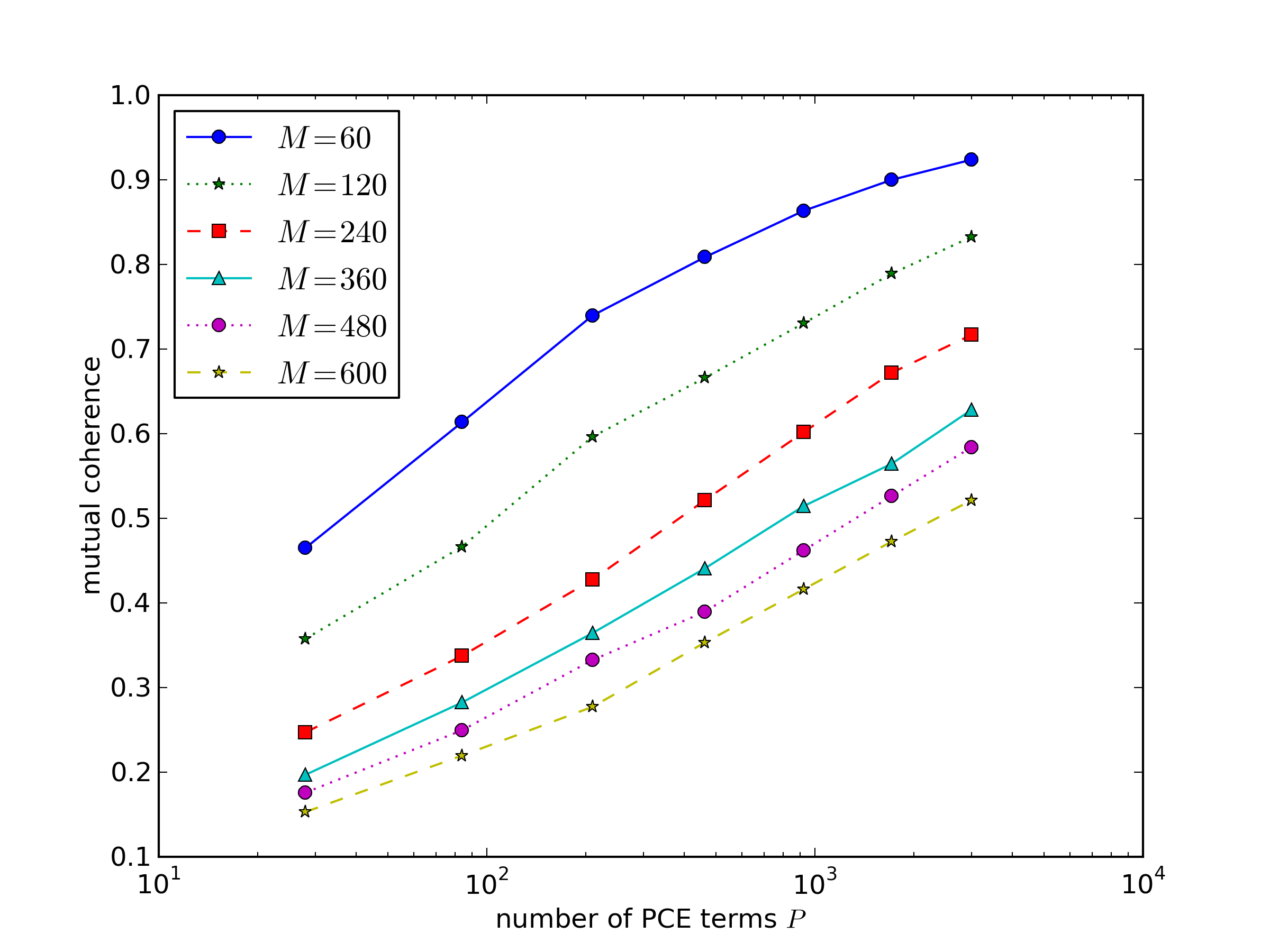}}\hspace{-.95cm}
\subfloat[]{\includegraphics[width=0.7\textwidth]{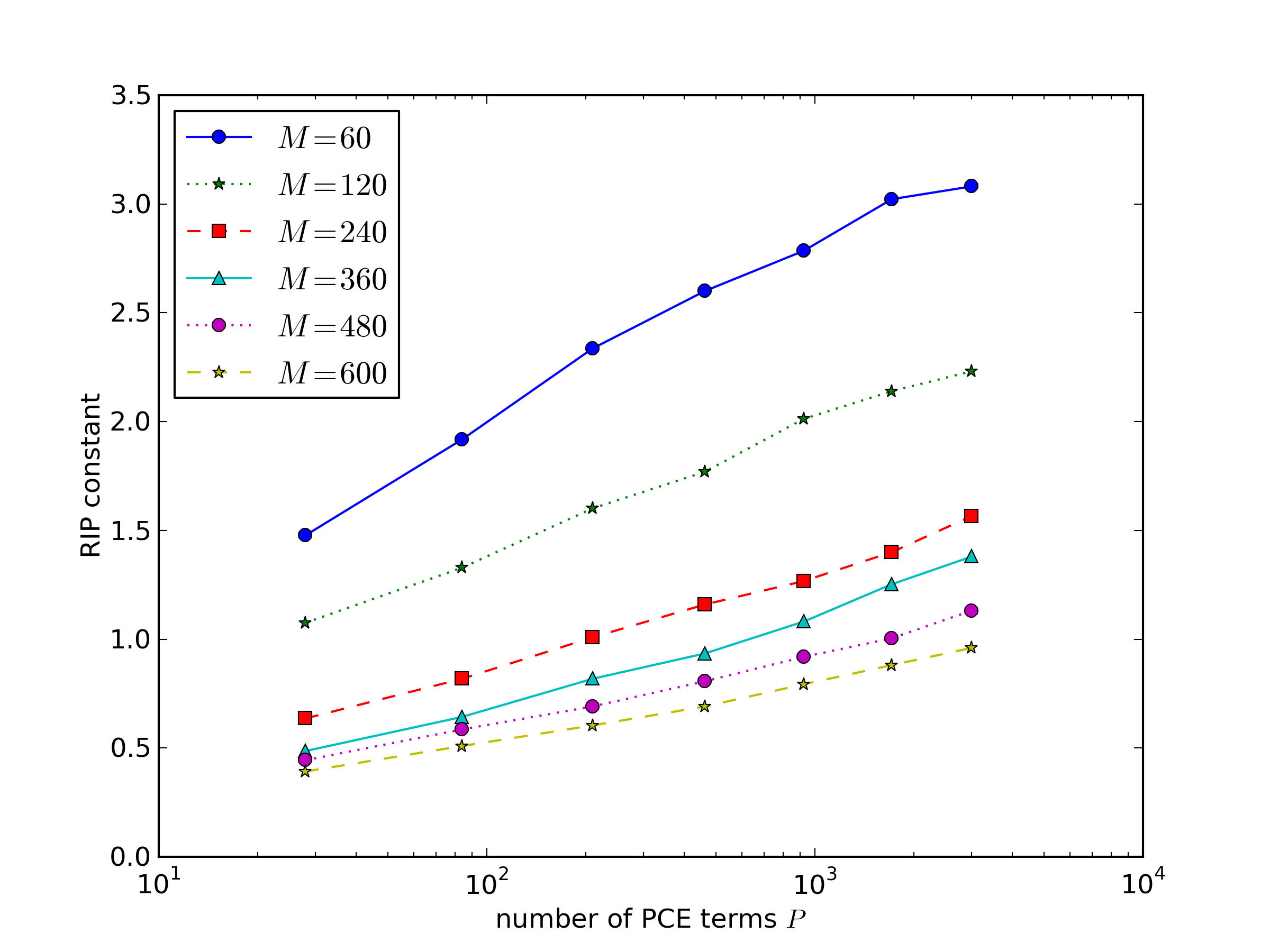}}
}
\caption{Growth of mutual coherence and RIP constant $\delta_{10}$ with the number of terms $P$ in a 6-dimensional total-degree multivariate Legendre PCE basis.}
\label{fig:mutual-coherence-increase}
\end{figure}

Naively choosing a large degree $p$ can cause a degradation in the accuracy of the PCE coefficients.
%Bounds on the number of samples $M$ needed to recover a Legendre PCE of a certain sparsity $s$ are given in~\cite{doostan2011}. 
%These bounds are not tight and known to be overly pessimistic, but what they stress is that when using a total degree or hyperbolic truncation the mutual coherence will increase
%as the degree $p$ is increased.
Figure~\ref{fig:mutual-coherence-increase} demonstrates that both the mutual coherence and the 10-sparse RIP constant $\delta_{10}$ of the Vandermonde matrix $\vand$ increases
as the number of basis terms $P$ increases.\footnote{The RIP constant reported here is a lower bound found by computing the eigenvalues of 10,000 randomly selected submatrices $\vand_{\bset_s}$.}

\begin{figure}[ht]
\centering
\includegraphics[width=0.95\textwidth]{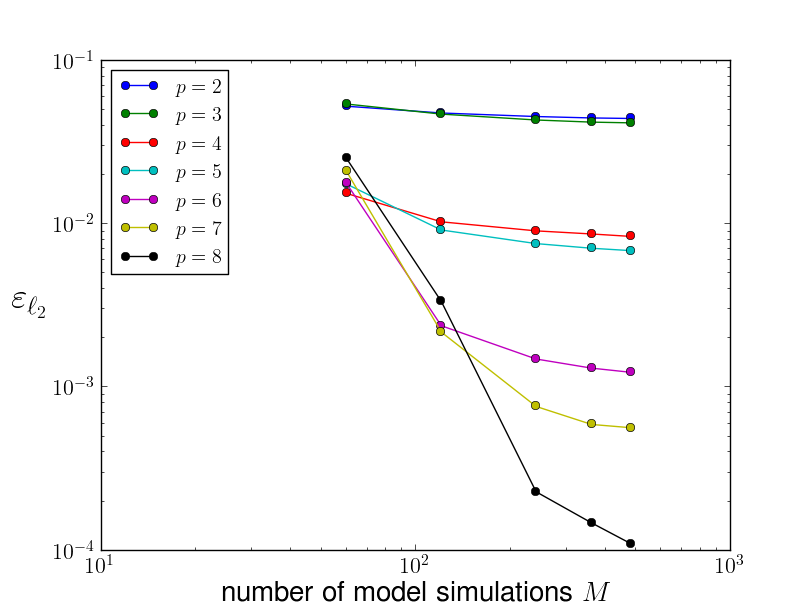}
%\includegraphics[width=0.45\textwidth]{transition-study-genz-os-no-decay-d-10-n-100.png}
%\caption{(Left) The dependence of PCE accuracy on the polynomial degree $p$ (Right) The degradation of recovery accuracy as the number of PCE terms increases.}
\caption{The dependence of PCE accuracy on the polynomial degree $p$. The $degree$ of the most accurate expansion is dependent on the number of LHS samples $M$ used to construct the PCE. Results 
were obtained using orthogonal matching pursuit with cross validation to choose $\stoptol$ applied to the 6-dimensional random oscillator~\eqref{eq:random-oscillator}. 
Accuracy is measured using the $\ell_2$ norm computed at $100,000$ Latin-hypercube samples (see Section~\ref{sec:results}).}
\label{fig:degree-study}
\end{figure}

Increases in mutual coherence and RIP constant are correlated with a decrease in the accuracy of PCE coefficients recovered by $\ell_1$-minimization.
Figure~\ref{fig:degree-study} demonstrates that, for a fixed number of samples, as the number of terms $N$ and, consequently, the mutual coherence and RIP constant increase 
(see Figure~\ref{fig:mutual-coherence-increase}), the PCE recovered by $\ell_1$-minimization becomes less accurate. Figure~\ref{fig:degree-study} also illustrates that the accuracy of the PCE depends upon the degree of the basis used. As the number of samples $M$
increases, $\ell_1$-minimization is able to recover more dominant coefficients and a higher degree should be used. 
However for a given number of samples increasing the degree does not always lead to a reduction in error. For example when $M=60$ the PCE of total degree $p=4$ has the smallest error and at $M=120$
setting $p=7$ produces the smallest error. It is not until $M=240$ that the highest degree basis $p=8$ produces the smallest error.
These results are consistent with the theoretical results in~\cite{doostan2011} that assert that 
the number of samples $M$ needed to recover a Legendre PCE of a certain sparsity $s$ increases with the number of terms in the PCE basis.

\section{Iterative basis selection}
\label{sec:iterative-basis-selection}
When the coefficients of a PCE can be well approximated by a sparse vector, $\ell_1$-minimization is extremely effective at recovering the coefficients of that PCE. 
It is possible, however, to further increase the efficacy of $\ell_1$-minimization by leveraging realistic models of structural dependencies between the values and 
locations of the PCE coefficients. For example~\cite{Baraniuk_CDH_IEEIT_2010,Duarte_WB_SPARS_2005,La_D_IEEEIP_2006} have successfully increased the performance
of $\ell_1$-minimization when recovering wavelet coefficients that exhibit a tree-like structure. In this vein, we propose an algorithm for identifying the large coefficients
of PC expansions that form a semi-connected subtree of the PCE coefficient tree.

\begin{figure}
\centering
\includegraphics[width=0.75\textwidth]{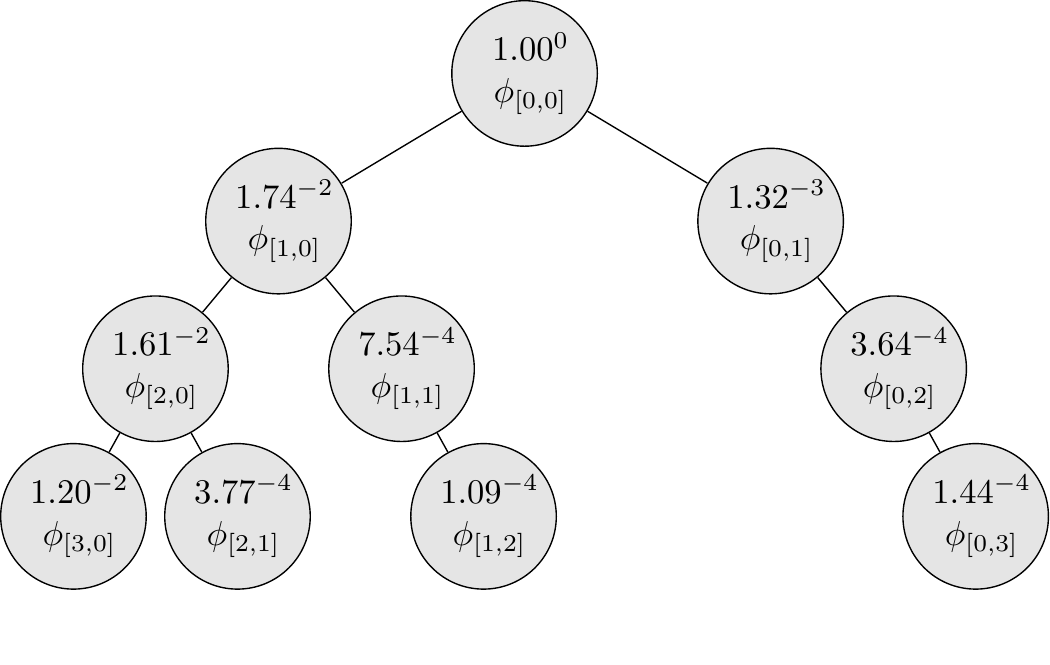}
\caption{Tree structure of the coefficients of a two dimensional PCE with $\bset_{3,1}^2$. For clarity we only depict one connection per node, but in $d$ dimensions a node of a given
degree $p$ will be a child of up to $d$ nodes of degree $p-1$. For example, not only is the basis $\bphi_{[1,1]}$ a child of $\bphi_{[1,0]}$ (as depicted) but it is also a child of 
$\bphi_{[0,1]}$} 
\label{fig:pce-tree}
\end{figure}

The coefficients of polynomial chaos expansions often form a multi-dimensional tree.
Given an ancestor basis term $\phi_{\blambda}$ of degree $\norm{\blambda}{1}$ we define the indices of its children as $\blambda+\be_k$, $k=1,\ldots,d$,
where $\be_k=(0,\ldots,1,\ldots,0)$ is the unit vector co-directional with the $k$-th dimension.
% We refer to the basis terms with $\hat{\blambda}-\be_k$ as ancestors of the basis indexed by $\hat{\blambda}$.
An example of a typical PCE tree is depicted in Figure~\ref{fig:pce-tree}. In this figure, as often in practice, the magnitude of the ancestors of a PCE coefficient is a
reasonable indicator of the size of the child coefficient. In practice, some branches (connections) between levels of the tree may be missing. We refer to trees with missing branches 
as semi-connected trees.

In the following we present a method for estimating PCE coefficients that leverages the tree structure of 
PCE coefficients to increase the accuracy of coefficient estimates obtained by $\ell_1$-minimization.

% \ks{I believe the tree structure is often seen in physics-driven models. One can always come up with synthetic cases that break this, but usually in physical models the coefficients decay in a certain way. I think also Alireza Doostan had some works where he showed coefficient decay in PDE solutions. Or that was only for 1d? Just my two cents.}

% Neither the total degree or hyperbolic truncation are ideal. The total degree truncation grows to quickly with dimension limiting its ability to recover high-degree terms.
% The hyperbolic index set will have difficulty recovering basis terms with many active variables. Moreover neither method is suited to the regime under which 
% the importance of the random variables varies anisotropically, that is higher degrees are needed in some dimensions than others. 
% To overcome these limitations we propose an iterative algorithm for choosing the PCE index set $\bset$.

\subsection{Algorithm}\label{sec:algorithm}
Typically $\ell_1$-minimization is applied to an a priori chosen and fixed basis set $\bset$. However the accuracy of coefficients obtained by $\ell_1$-minimization can be increased by
adaptively selecting the PCE basis.

\begin{figure}
\hspace{-2cm}\includegraphics[width=1.3\textwidth]{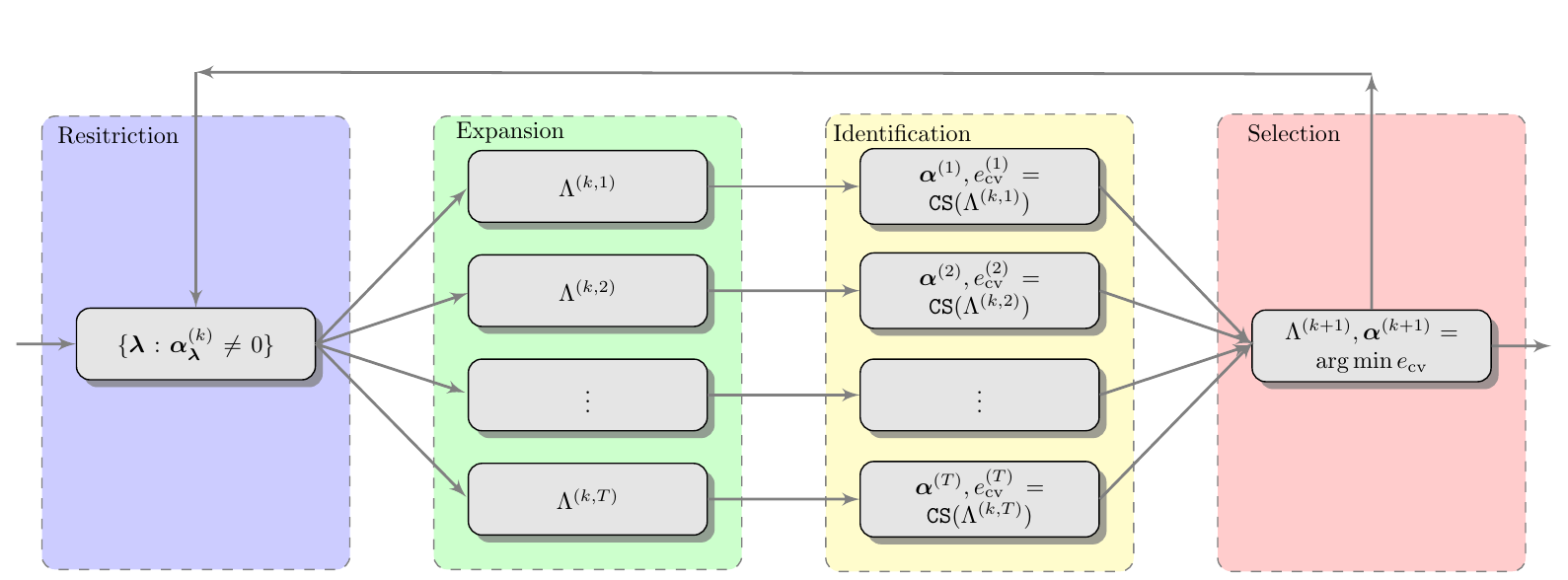}
\caption{Graphical depiction of the basis adaptation algorithm.}
\label{fig:basis-selection-alg}
\end{figure}

To select a basis for $\ell_1$-minimization we employ a four step iterative procedure involving restriction, expansion, identification and selection. 
% The algorithm begins with 
% an initial basis and corresponding coefficients. This basis is then restricted to only contain basis terms with non-zero coefficients. This restricted basis is
% expanded \ks{say how} to generate a set of new candidate bases. $\ell_1$-minimization ($\ell_1$-minimization) and cross validation is then used to identify a set of large coefficients 
% that provide a good balance between fitting and prediction \ks{for each of the candidate basis set?}. In the last step cross validation error is used to select the best candidate basis. These steps
% are repeated until the new candidate bases generated no longer reduce the cross validation error. \ks{Frankly, this verbal explanation was not too clear. Maybe identify after each sentence in parentheses which of the four steps you refer to as it goes. Of course Figure 6 clarifies it though:)}
%Hence forth we will refer to the three aforementioned steps as selection, restriction and expansion. 
The iterative basis selection procedure is outlined in Algorithm~\ref{alg:basis-selection}. A graphical version of the algorithm is also presented in Figure~\ref{fig:basis-selection-alg}.
The latter emphasizes the four stages of basis selection, that is restriction, growth, identification and selection. These four stages are also highlighted in 
Algorithm~\ref{alg:basis-selection} using the corresponding colors in Figure~\ref{fig:basis-selection-alg}.

To initiate the basis selection algorithm, we first define a basis set $\Lambda^{(0)}$ and use $\ell_1$-minimization to identify the largest coefficients $\coef^{(0)}$. The choice of $\Lambda^{(0)}$
can sometimes affect the performance of the basis selection algorithm. We found a good choice to be $\Lambda^{(0)}=\bset_{p,1}$,
where $p$ is the degree that gives $\card{\bset^d_{p,1}}$ closest to $10M$, i.e. $\bset^d_{p,1} = \argmin_{\bset^d_{p,1}\in\{\bset^d_{1,1},\bset^d_{2,1},\ldots\}}\abs{\card{\bset^d_{p,1}}-10M}$.
Given a basis $\Lambda^{(k)}$ and corresponding coefficients $\coef^{(k)}$ we reduce the basis to a set $\Lambda^{(k)}_\stoptol$ containing only the terms with non-zero coefficients. 
This restricted basis is then expanded $T$ times using an algorithm which we will describe in Section~\ref{sec:basisexp}. $\ell_1$-minimization is then applied to each of the expanded basis 
sets $\Lambda^{(k,t)}$ for $t=1,\dots, T$.
Each time $\ell_1$-minimization is used, we employ cross validation to choose $\stoptol$. Therefore, at every basis set considered during the evolution of the algorithm we have a measure
of the expected accuracy of the PCE coefficients. At each step in the algorithm we choose the basis set that results in the lowest cross validation~error.

% \begin{algorithm}[ht]
% \caption{$\Lambda^\star$,$\coef^\star$=BASIS\_SELECTION[$\vand$,$\rhs$,$\stoptol$] }
% Set $\Lambda^{\star} = \Lambda_0 = \bset^d_{p,1} = \argmin_{\bset^d_{p,1}\in\{\bset^d_{1,1},\bset^d_{2,1},\ldots\}}\abs{\card{\bset^d_{p,1}}-10M}$\\
% Set  $t^\star=3$, $\cverr^\star = \infty$, $i = 1$, \\
% While TRUE
% \begin{itemize}
% \item $\coef_i$, $\cverr_i$ = $\ell_1$-minimization[$\vand(\Lambda_{i-1})$,$\rhs$]
% \item $\Lambda_i^\stoptol=\{\blambda:\blambda\in\Lambda_{i-1}, \coef_{\blambda} \ne 0\}$
% \item for $t\in\{1,\ldots,t^\star\}$
% \begin{itemize}
% \item $\Lambda_{i,t}$ = EXPAND[$\Lambda_i^\stoptol(\coef_i)$,$t$]
% \item $\coef_{i,t}$, $\cverr_{i,t}$ = $\ell_1$-minimization[$\vand(\Lambda_{i,t})$,$\rhs$]
% \item if $\cverr_{i,t}<\cverr_i$: $\cverr_i=\cverr_{i,t}$, $\coef_i=\coef_{i,t}$, $\Lambda_i=\Lambda_{i,t}$, $t^\star=t$
% \end{itemize}
% 
% \item if $\cverr_i^\star > \cverr^\star$ terminate
% \item $\cverr^\star=\cverr_i$, $\Lambda^\star=\{\blambda:\blambda\in\Lambda_{i}, \coef_{\blambda} \ne 0\}$, $\coef^\star = \{\coef_{\blambda} : \blambda\in\Lambda_{i} ,\coef_{\blambda}\ne 0\}$
% \item $i=i+1$
% \end{itemize}
% \label{alg:basis-selection}
% \end{algorithm}

\begin{algorithm}[H]
%\dontprintsemicolon%
\DontPrintSemicolon % Some LaTeX compilers require you to use \dontprintsemicolon instead
\footnotesize
$\Lambda^{\star} = \Lambda^{(0)} = \bset^d_{p,1} = \argmin_{\bset^d_{p,1}\in\{\bset^d_{1,1},\bset^d_{2,1},\ldots\}}\abs{\card{\bset^d_{p,1}}-10M}$\;
$\coef^{(0)}$, $\cverr^{(0)}$ = $\ell_1$-minimization[$\vand(\Lambda^{(0)})$,$\rhs$]\;
$T=3$, $\cverr^\star = \infty$, $k = 1$\;
\While{TRUE}{
  $\cverr^{(k)}=\infty$\;
  \RHiLi$\Lambda^{(k,0)}=\{\blambda:\blambda\in\Lambda^{(k-1)}, \coef_{\blambda}^{(k)} \ne 0\}$\;
  \For{$t\in\{1,\ldots,T\}$}{
     \EHiLi$\Lambda^{(k,t)}$ = EXPAND[$\Lambda^{(k,t-1)}$]\;
     \IHiLi$\coef^{(k,t)}$, $\cverr^{(k,t)}$ = $\ell_1$-minimization[$\vand(\Lambda^{(k,t)})$,$\rhs$]\;
     \If {$\cverr^{(k,t)}<\cverr^{(k)}$} {\SHiLi$\cverr^{(k)}=\cverr^{(k,t)}$, $\coef^{(k)}=\coef^{(k,t)}$, $\Lambda^{(k)}=\Lambda^{(k,t)}$}
 }
  \If{$\cverr^{(k)} > \cverr^\star$}{ TERMINATE }
  $\alpha^\star=\coef^{(k)},\; \Lambda^\star=\Lambda^{(k)},\; \cverr^\star=\cverr^{(k)}$\;%, $\Lambda^\star=\{\blambda:\blambda\in\Lambda^{(k)}, \coef_{\blambda} \ne 0\}$, $\coef^\star = \{\coef_{\blambda} : \blambda\in\Lambda^{(k)} ,\coef_{\blambda}\ne 0\}$\;
%   $k=k+1$\;
}
\caption{$\Lambda^\star$,$\coef^\star$=BASIS\_SELECTION[$\vand$,$\rhs$,$\stoptol$] }
\label{alg:basis-selection}
\end{algorithm}

\subsubsection{Basis expansion} \label{sec:basisexp}

Define $\{\blambda+\be_j:1\le j\le d\}$ the forward neighborhood of an index $\blambda$ and similarly let $\{\blambda-\be_j:1\le j\le d\}$ denote the backward neighborhood.
To expand a basis set $\Lambda$ we must first find the forward neighbors $\mathcal{F}=\{\blambda+\be_j : \blambda\in\Lambda, 1\le j\le d \}$ of all indices $\blambda\in\Lambda$.
The expanded basis is then given by 
\[
\Lambda^+=\Lambda\cup\mathcal{A},\quad \mathcal{A}=\{\blambda: \blambda\in\mathcal{F}, \blambda-\mathbf{e}_n\in\Lambda\text{ for }1\le n\le d,\, \lambda_k > 1\}
\]
where we have used the following admissibility criteria 
\begin{equation}
\label{eq:admissibility}
\blambda-\mathbf{e}_n\in\Lambda\text{ for }1\le n\le d,\, \lambda_k > 1
\end{equation}
to target PCE basis indices that are likely to have large PCE coefficients. A forward neighbor is admissible only if its backward neighbors exist in all dimensions. 
If the backward neighbors do not exist then $\ell_1$-minimization has previously identified that the coefficients of these backward neighbors are negligible. 

\begin{figure}[ht]
\begin{center}
\includegraphics[width=0.95\textwidth]{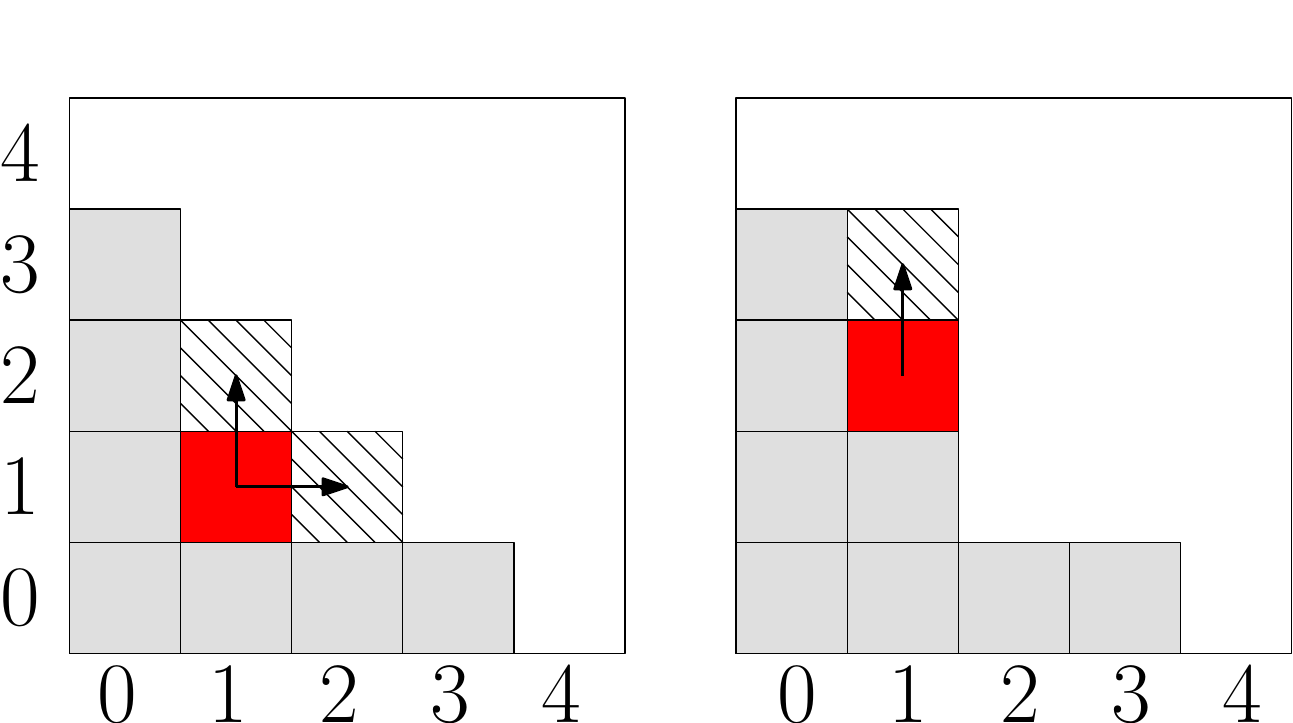}
\caption{Identification of the admissible indices of an index (red). The indices of the current basis $\Lambda$ are gray and admissible indices are striped. A index is admissible only if
its backwards neighbors exists in every dimension.}
\label{fig:index-dmissibiliy-examples}
\end{center}
\end{figure}

The admissibility criterion is explained graphically in Figure~\ref{fig:index-dmissibiliy-examples}. In the left graphic, 
both children of the current index are admissible, because its backwards neighbors exist in every dimension. 
In the right graphic only the child in the vertical dimension is admissible,
as not all parents of the horizontal child exist. 
% \begin{algorithm}[H]
% \DontPrintSemicolon % Some LaTeX compilers require you to use \dontprintsemicolon instead
% \footnotesize
% \For{$\blamda\in\Lambda$}{
% $\mathcal{F}=\{\blambda+\be_j : \blambda\in\Lambda, 1\le j\le d \}$
% $\mathcal{A}=\{\blambda: \blambda\in\mathcal{F}, \blambda-\mathbf{e}_n\in\Lambda\text{ for }1\le n\le d,\, \lambda_k > 1\}$
% $\Lambda^+=\Lambda\cup\mathcal{A}$
% }
% \caption{$\Lambda^{+}$ = EXPAND[$\Lambda$] }
% \label{alg:basis-selection}
% \end{algorithm}

At the $k$-th iteration of Algorithm~\ref{alg:basis-selection}, $\ell_1$-minimization is applied to $\Lambda^{(k-1)}$ and used to identify the significant coefficients of the PCE and their corresponding basis terms $\Lambda^{(k,0)}$. 
The set of non-zero coefficients $\Lambda^{(k,0)}$ identified by $\ell_1$-minimization is then expanded. 
% The key to the proposed method working well is to ensure the basis generated by expansion is able capture higher degree terms without increasing the mutual coherence
% to a point which degrades the ability of $\ell_1$-minimization to recover those higher order terms. 
The \texttt{EXPAND} routine expands an index set by one polynomial degree,
but sometimes it may be necessary to expand the basis $\Lambda^{(k)}$ more than once.\footnote{The choice of $T>1$ enables the basis selection algorithm to be applied to semi-connected 
tree structures as well as fully connected trees. Setting $T>1$ allows us to prevent premature termination of
the algorithm if most of the coefficients of the children of the current set $\bset^{(k)}$ are small but the coefficients of the children's children are not.}
%For example by more than one degree to avoid situations when no basis terms one degree higher are significant, but basis terms two or three degrees higher are. 
To generate these higher degree index sets \texttt{EXPAND} is applied recursively to $\Lambda^{(k,0)}$ up to a fixed number of $T$ times.
Specifically, the following sets are generated
$$\Lambda^{(k,t)}=\Lambda^{(k,t-1)}\cup\{\blambda:\blambda-\mathbf{e}_n\in\Lambda^{(k,t-1)},1\le n\le d,\, \lambda_n > 1\}.$$  
As the number of expansion steps $T$ increases the number of terms in the expanded basis increases rapidly and
degradation in the performance of $\ell_1$-minimization can result (this is similar to what happens when increasing the degree of a total degree basis). 
To avoid degradation of the solution, we use cross validation to choose the number of inner expansion steps $t\in[1,T]$.

\section{Numerical examples}\label{sec:results}
In this section we use several numerical tests to demonstrate the benefit of the basis selection method.
In each example we seek a PCE approximation to a model output given a
set of uncertain parameters with a known range or distribution. 
% Without loss of generality we consider variables with 
% a finite range and use one-dimensional Legendre Polynomials to build the final multivariate 
% polynomial approximations. The method proposed can easily be extended to variables with 
% infinite range using appropriate one-dimensional polynomials, following the typical 
% Askey scheme~\cite{xiu02a}.

We compare the approximations constructed
using basis selection against those constructed using a non-adaptive strategy. 
The non-adaptive strategy consists of generating basis sets $\bset_{2,1}^d,\ldots,\bset_{p,1}^d$ where
$p$ is the degree that produces the basis set with a cardinality closest to $100000$. $\ell_1$-minimization 
is then applied to this basis with a cross validation tolerance search to compute the non-zero polynomial
coefficients. The resulting basis with the lowest cross validation error is chosen to be the final approximation.

We also compare the non-adaptive and basis selection methods against OMP using a basis oracle. For
a set of $M$ model runs, the oracle method sets $\bset$ in~\eqref{eq:bpdn} to be the basis of the best $M$-term PCE approximation of the function.
The best $M$-term approximation is obtained by using a dimension adaptive sparse grid to calculate the `exact' PCE
coefficients and selecting the basis terms with the $M$ largest coefficients. This basis will be close to optimal
and therefore will serve as a good estimate of the maximum accuracy that can be gained from the use of basis selection.

In order to construct a PCE approximation, the sample design has to be specified. 
By a design we mean the choice of sample size, $M$, and the
selection $\bXi=\{\brv_i\}_{i=1}^M$
Here we opt for uniform random samples of size $M$. For small sample sizes $M$ 
the selection of $\bXi$ significantly affects the performance
of any approximation method. Therefore, for each $M$, twenty different designs are used to build
a PCE and the mean, maximum and minimum of the resulting errors
are reported. We did investigate the utility of using samples drawn from 
the Chebyshev measure but found that there was no consistent benefit. Even in the cases for which a benefit was observed, 
the improvement was small relative to the benefit gained from using basis adaptation.
This finding is consistent with~\cite{Yan_GX_IJUQ_2012}.

To measure the performance of an approximation, we will use the $\ell_2$ error (RMSE).
% and $l_\infty$ (Max) error of the approximations.
Specifically given a set of $Q=100,000$ Latin-hypercube samples $\bXi_{\text{test}}=\{\brv^{(i)}\}_{i=1}^Q\in\dom$ and samples of the
true function $f(\brv^{(i)})$ and the PCE approximation $\hat{f}(\brv^{(i)})$ we compute
\[
 \varepsilon_{\ell_2} = \left(\frac{1}{Q}\sum_{i=1}^Q \lvert\hat{f}(\brv^{(i)})-f(\brv^{(i)})\rvert^2\right)^{1/2}
\]
Note in all examples presented using Legendre polynomials we transform each $d$ dimensional parameter domain $\dom$, build points $\bXi$ and test points $\bXi_{\mathrm{test}}$ to $[-1,1]^d$.
% Surrogates are also often used to estimate statistical moments such as mean and variance, 
% with respect to the uncertain inputs. Consequently we will also measure the performance of a 
% surrogate with the two additional measures
% \[
%  \varepsilon_{\mu} = \frac{\mean{\hat{f}} - \mean{f}}{\mean{f}}\quad\text{and}\quad 
% \varepsilon_{\sigma^2}  = \frac{\var{\hat{f}} - \var{f}}{\var{f}}
% \]
% where the mean and variance are computed using $\bXi_{\text{test}}$.

\subsection{Algebraic test function}\label{sec:genz-function-based-cs}
Consider the algebraic corner-peak test function~\cite{genz87}
\begin{equation}\label{eq:genz-corner-peak}
 f_{\mathrm{CP}}(\bx)=\left(1+\sum_{k=1}^d c_k\, \rv_k \right)^{-(d+1)},\quad \brv\in[0,1]^d
\end{equation}
This function provides a flexible test that can be used to identify the strengths of the proposed
algorithm. Specifically, the coefficients $c_k$ can be used to control the effective dimensionality 
and the compressibility of these functions. Here we will examine performance using three different
choices of $\bc = (\subs{c}{1}{d})^T$, specifically 

\[
c^{(1)}_k=\frac{k-\frac{1}{2}}{d},\quad c_k^{(2)}=\frac{1}{k^2}\quad \text{and} 
\quad c_k^{(3)} = \exp\left(\frac{k\log(10^{-8})}{d}\right), \quad k=1,\ldots,d
\]
normalizing such that $\sum_{k=1}^d c_k = 0.25$. The coefficients $\bc^{(1)}$, $\bc^{(2)}$ and $\bc^{(3)}$
represent increasing levels of anisotropy and decreasing effective dimensionality. Anisotropy refers to the dependence of the function variability,
often measured through variance, on individual parameter dimensions $\rv_n$. When a function is strongly anisotropic, the majority of the function variance
can be attributed to a small set of dimensions. The size of this subset is referred to as the effective dimension.

% Before demonstrating the efficacy of the proposed basis selection method it is important to show the importance of choosing $\stoptol$ well.
% Figure~\ref{fig:truncation-tolerance-cross-validation-vs-l2-error}
% presents a typical example illustrating the change in the $\ell_2$ error of a PCE with fixed degree, as the tolerance $\stoptol$ is decreased. 
% The figure also plots the cross validation error which is a good estimation of the $\ell_2$ error. The vertical line represents the tolerance chosen
% by cross validation and the horizontal line is the $\varepsilon_{\ell_2}$ error in the resulting PCE. The result shown is typical. 
% There is a bias in the cross validation estimate, yet despite this bias cross validation consistently chooses a tolerance that produces a near minimal error.
% \begin{figure}
% \centering
% \includegraphics[width=0.95\textwidth]{tolerance-cv-comparison-genz-cp-no-decay-d-10-n-200-p-4.png}
% \caption{The use of cross validation to select the truncation tolerance $\stoptol$ for $c^{(1)}$. The vertical line represents the tolerance chosen
% by cross validation and the horizontal line is the $\varepsilon_{\ell_2}$ error in the resulting PCE. $M=200$ samples were used.}
% \label{fig:truncation-tolerance-cross-validation-vs-l2-error}
% \end{figure}
 
\begin{figure}
% the use of makebox allows figures to spill into margins and be centered despite spill
\makebox[1\linewidth][c]{%
\subfloat[]{\includegraphics[width=0.7\textwidth]{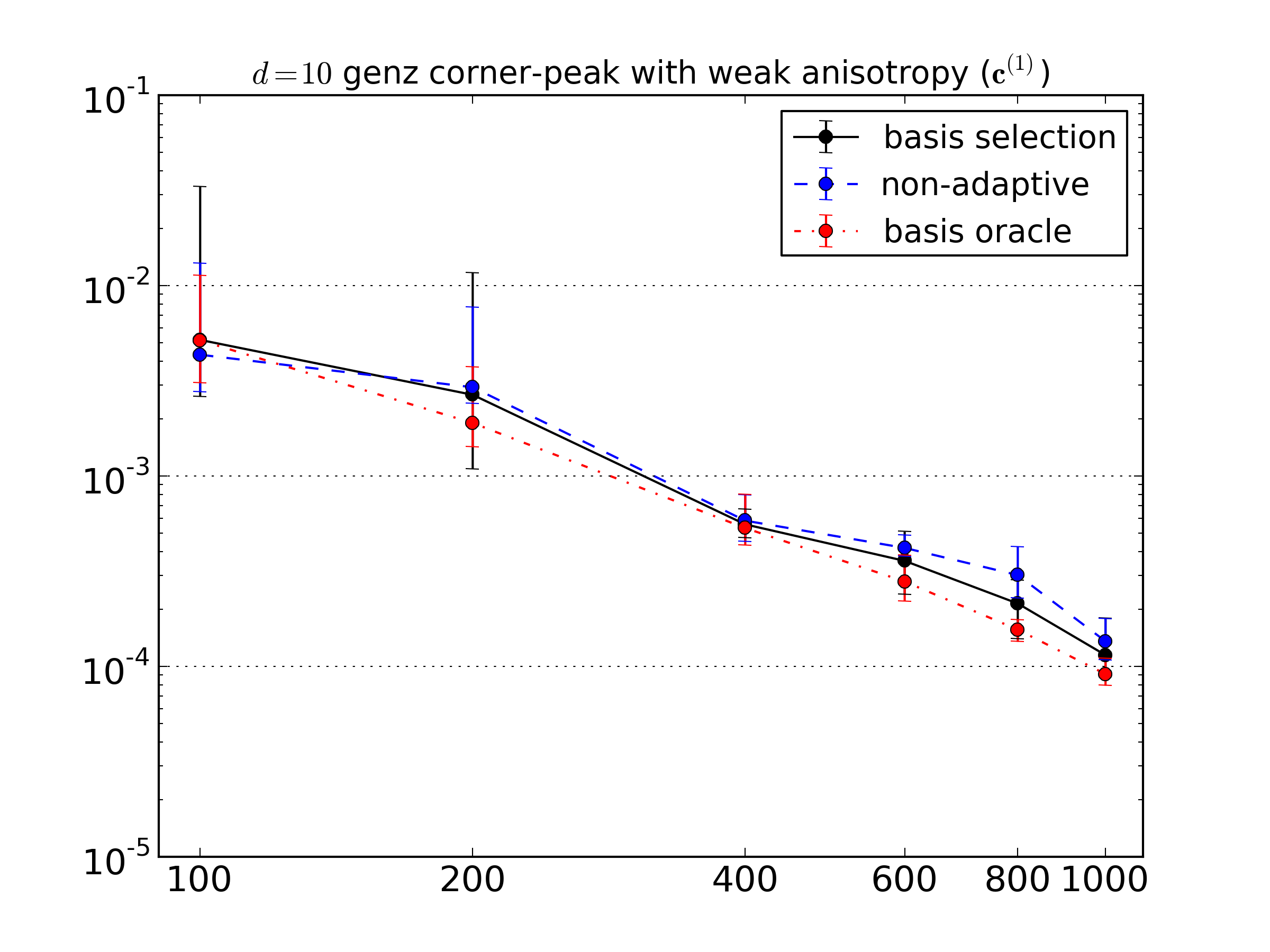}}
\hspace{-.95cm}
\subfloat[]{\includegraphics[width=0.7\textwidth]{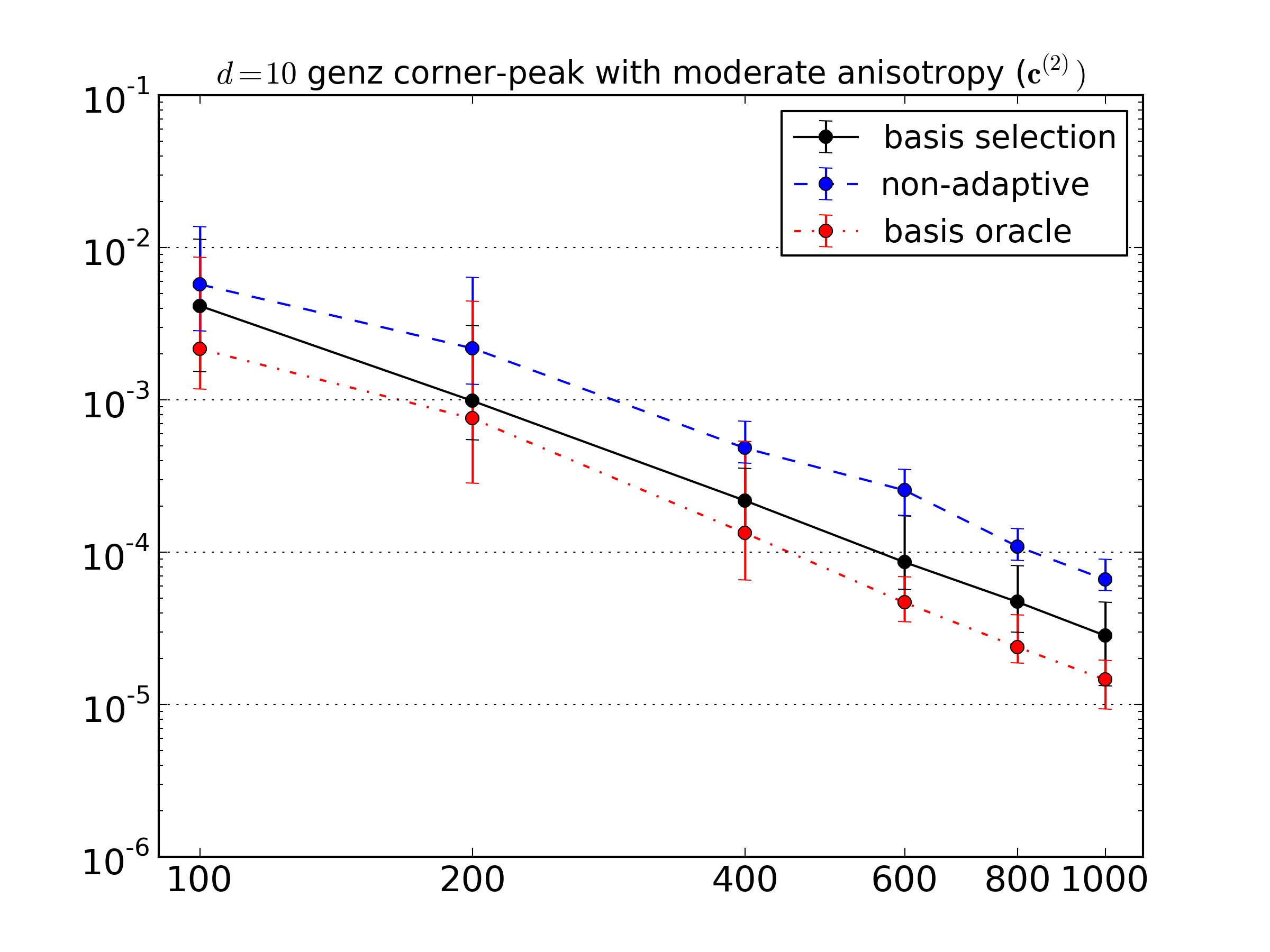}}
}
\makebox[1\linewidth][c]{%
\subfloat[]{\includegraphics[width=0.7\textwidth]{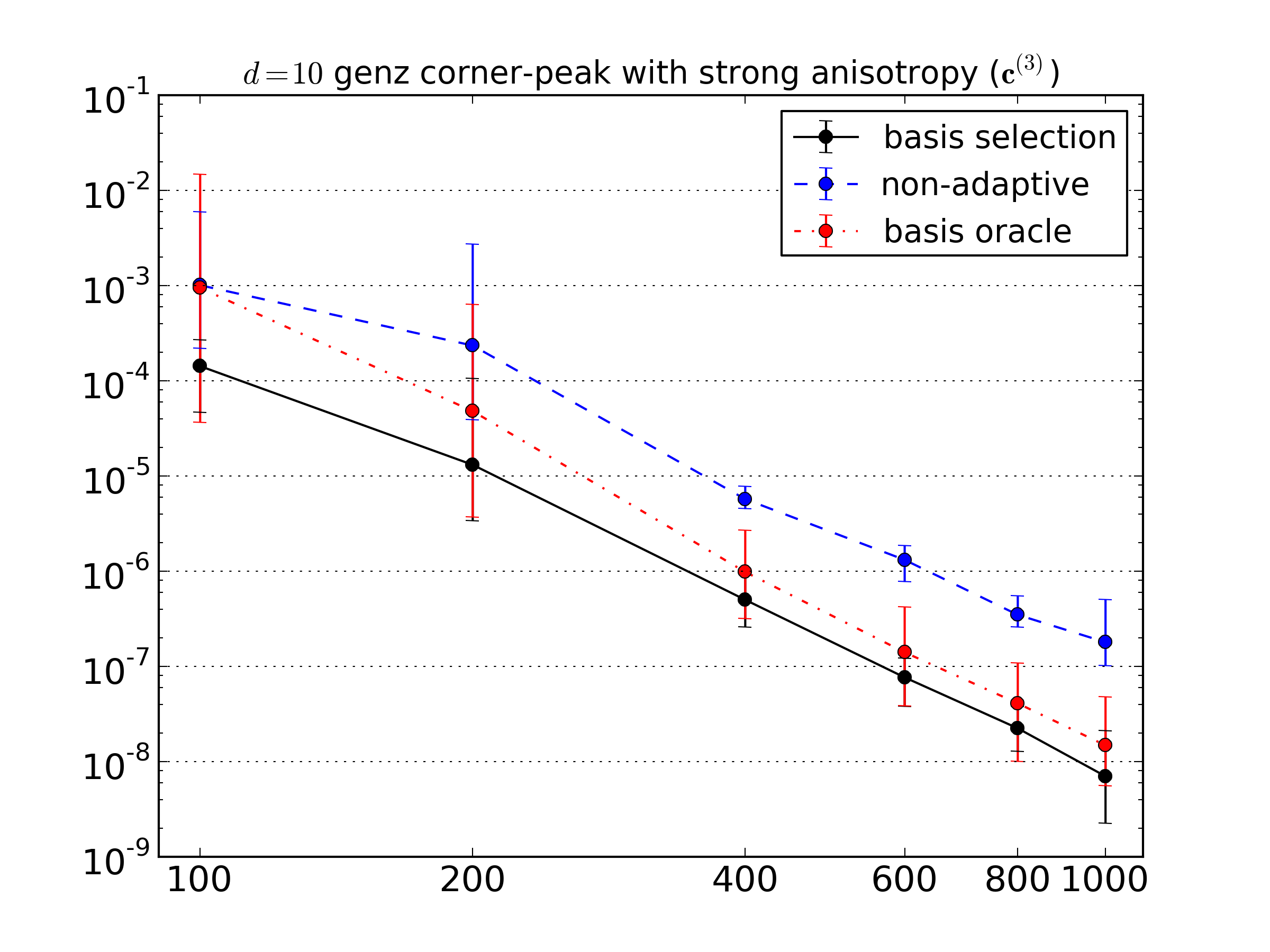}}\hspace{-.95cm}
\subfloat[]{\includegraphics[width=.7\textwidth]{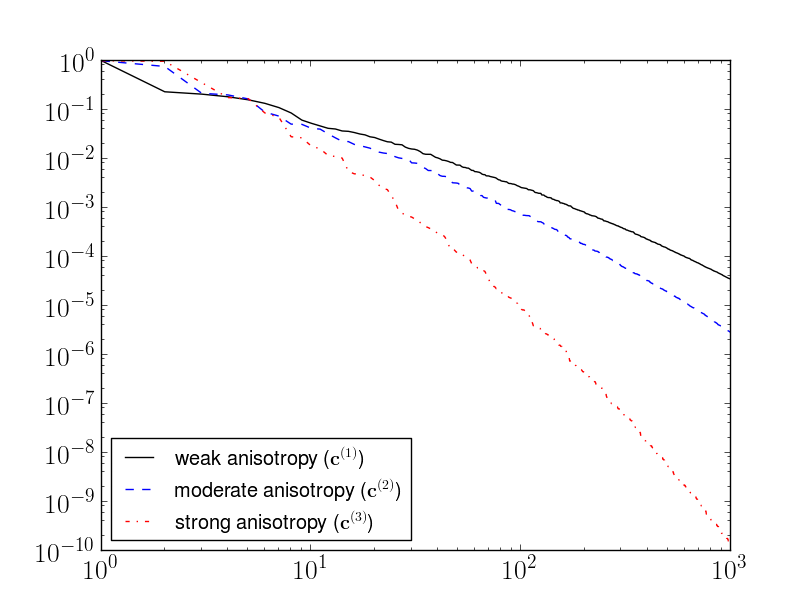}}
}
\caption{Convergence of the RMSE, with respect to increasing design size $M$, in the Legendre PCE approximation of the Genz corner-peak 
function~\eqref{eq:genz-corner-peak} for the three coefficient regimes: 
(a) $\bc^{(1)}$, (b) $\bc^{(2)}$, (c) $\bc^{(3)}$. 
The error bars represent the minimum and maximum error over the $20$ trials for each design size $M$.
(d) PCE coefficients, sorted by magnitude, for the three coefficient regimes.
}
\label{fig:genz-various-decays}
\end{figure}

The performance of the basis selection method is dependent on the properties of the model being approximated.
Figure~\ref{fig:genz-various-decays} plots the $\varepsilon_{\ell_2}$ error in
the polynomial approximations for increasing number of model evaluations.
For all three levels of anisotropy the adaptive method produces an expansion no worse than the non-adaptive method, for the same 
sample size. When anisotropy is introduced the accuracy of basis selection increases relative to the non-adaptive method. The stronger the anisotropy the better the 
relative performance. Figure~\ref{fig:genz-various-decays} also plots the PCE obtained using OMP with an oracle basis. When very weak anisotropy is present $\bc^{(1)}$, there is little 
that can be gained by using a well chosen basis, as evident by the lack of separation between the three convergence curves. However as the strength of the anisotropy is increased
the effect of the oracle basis on accuracy becomes much more apparent. When strong anisotropy, $\bc^{(3)}$, is present, basis selection is able to obtain the same
accuracy as the oracle without a priori information on the truncation of the basis which is required by the oracle. In the moderately anisotropic $\bc^{(2)}$ 
case basis selection does not perform as well as the oracle but does still perform
better than the non-adaptive method.

% \begin{figure}[!htb]
% \centering
% \includegraphics[width=0.9\textwidth]{genz-sorted-coeff-decay.png}
% \caption{PCE coefficients, sorted by magnitude, for the three coefficient regimes.}
% \label{fig:genz-sorted-coeff-decays}
% \end{figure}
To understand the correlation between anisotropy and the performance of the basis selection method we must consider the the structure of the PCE coefficients induced by varying $\bc$.
Figure~\ref{fig:genz-various-decays} (d) plots the decay of the PCE coefficients when sorted by magnitude. As anisotropy increases, so does the rate of decay of the sorted PCE coefficients.

It is the strength of decay that controls the performance of basis selection. When the rate of decay is high then the function is more compressible and thus better suited 
to being approximated using $\ell_1$-minimization. Anisotropy will often result in compressible coefficients, but it is conceptually possible for models to be compressible without being anisotropic.
For some problems such as the elliptic Poisson equation the decay of the PCE coefficients can be calculated a-priori~\cite{Babuska_TZ_SIAMNA_2004,Beck_TNT_MMMAS_2012,Cohen_DS_FCM_2010}
but unfortunately in practice, the decay of the PCE coefficients of a model cannot be determined ahead of time. A practical means of identifying the coefficient decay regime would be very useful 
but is beyond the scope of this paper.

Not only does the rate of coefficient decay affect performance, but so to does the ability of the PCE basis $\bset$ to represent the target function.
For example if the `true' PCE has large high degree terms with large coefficients but the basis $\bset$ does not have these high degree terms then the PCE obtained using
$\bset$ will not be as accurate as a PCE obtained using a basis that included the important high degree terms.

\begin{figure}[h!]
\centering
\makebox[1\linewidth][c]{%
\subfloat[]{\includegraphics[width=0.7\textwidth]{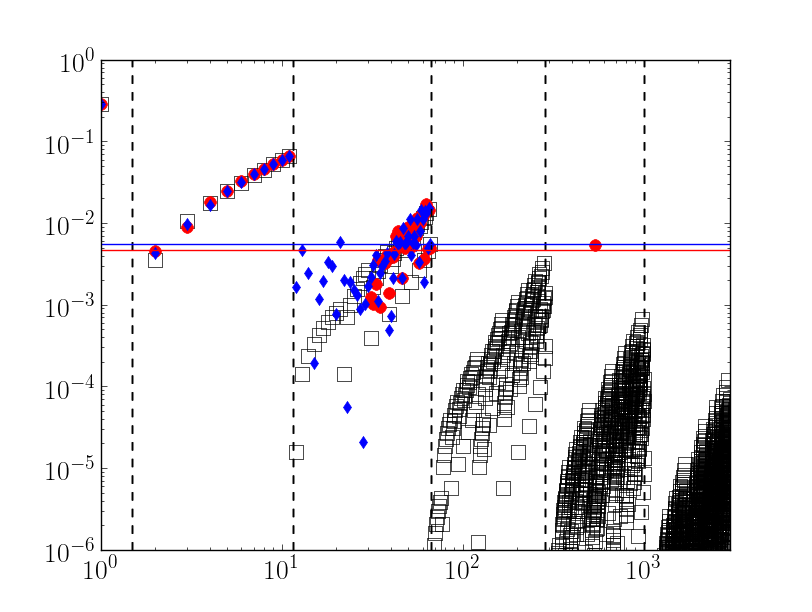}}
\hspace{-.95cm}
\subfloat[]{\includegraphics[width=0.7\textwidth]{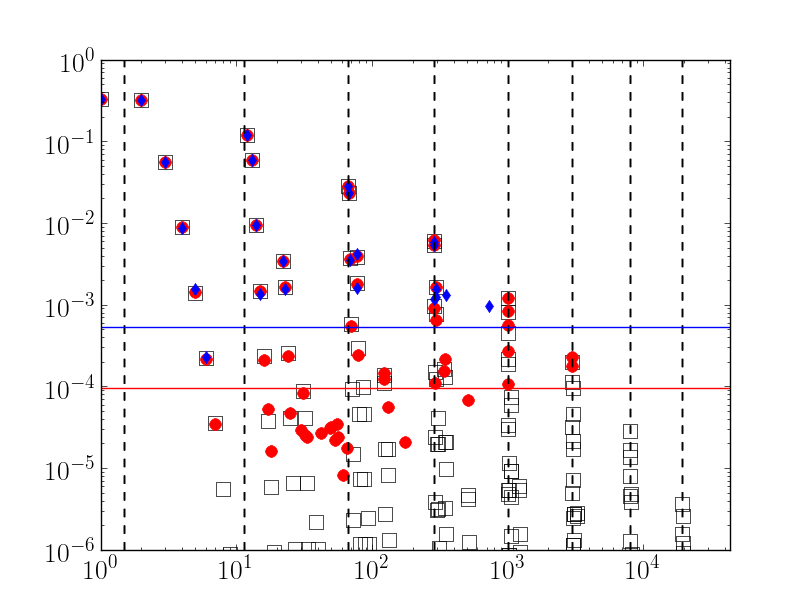}}
}
\caption{Comparison of PCE coefficients of $\bc^{(1)}$ (a) and $\bc^{(3)}$ (b). The black squares represent the `exact' coefficients, the red circles the non-zero coefficients of the 
basis selection method and the blue diamonds are the non-zero coefficients recovered by the non-adaptive approach. The indices on the x axis are sorted lexicographically by degree. The dashed vertical lines separate the PCE terms into degrees. The horizontal lines
represent the $\ell_2$ error in the basis selection and non-adapted PCE. The identification of more terms accurately with basis selection results in a smaller error.
}
\label{fig:pce-coefficient-regimes} 
\end{figure}

Figure~\ref{fig:pce-coefficient-regimes} plots the exact PCE coefficients of the algebraic test function using $\bc^{(1)}$ and $\bc^{(3)}$.
The `exact' coefficients obtained using a dimension-adaptive sparse grid with $100,000$ which results in an approximation error below $10^{-8}$.
In Figure~\ref{fig:pce-coefficient-regimes} (a), the random variables contribute similarly to the total variance of $f_\text{CP}$ 
and so the dominant coefficients are concentrated
in the lower degree terms of the PCE, thus a total degree basis set will perform as well as any alternative. 
In comparison, the importance of the dimensions of the function shown in Figure~\ref{fig:pce-coefficient-regimes} (b) decay exponentially with dimension, 
which results in higher-degree terms with large coefficients in some dimensions. In this coefficient regime, if $\ell_1$-minimization can only be applied with a low degree polynomial 
(which is true when using a total degree basis), the accuracy of the resulting PCE, for a given number of samples $M$, will not be as high as a PCE constructed using 
a basis that includes the dominant high-degree terms. Basis selection will typically allow identification and recovery of more high-degree coefficients than would be possible 
if using a total degree basis.

\subsection{Random oscillator}
In this section we investigate the performance of basis selection to quantify uncertainty in a damped 
linear oscillator subject to external forcing with six unknown 
parameters. That is,
\begin{equation}\label{eq:random-oscillator}
\frac{d^2x}{dt^2}(t,\brv)+\gamma\frac{dx}{dt}+k x=f\cos(\omega t),
\end{equation}
subject to the initial conditions
\begin{equation}
x(0)=x_0,\quad \dot{x}(0)=x_1,
\end{equation}
where we assume the damping coefficient $\gamma$, spring constant $k$,
forcing amplitude $f$ and frequency $\omega$, and the initial
conditions $x_0$ and $x_1$ are all uncertain. We solve~\eqref{eq:random-oscillator} analytically
to avoid consideration of discretization errors in our study.

Defining $
\brv=(\gamma,k,f,\omega,x_0,x_1)$ let $\rv_1\in[0.08,0.12]$, $\rv_2\in[0.03,0.04]$, $\rv_3\in[0.08,0.12]$, $\rv_4 \in [0.8,1.2]$, 
$\rv_5\in[0.45,0.55]$, $\rv_6\in[-0.05,0.05]$. For any parameter realization in $I_\brv$ the
harmonic oscillator will be underdamped. In the following, we choose our quantity of interest
to be the position $x(t)$ of the oscillator at $t=20$ seconds.

\begin{figure}
\centering
\makebox[1\linewidth][c]{%
\subfloat[]{\includegraphics[width=0.7\textwidth]{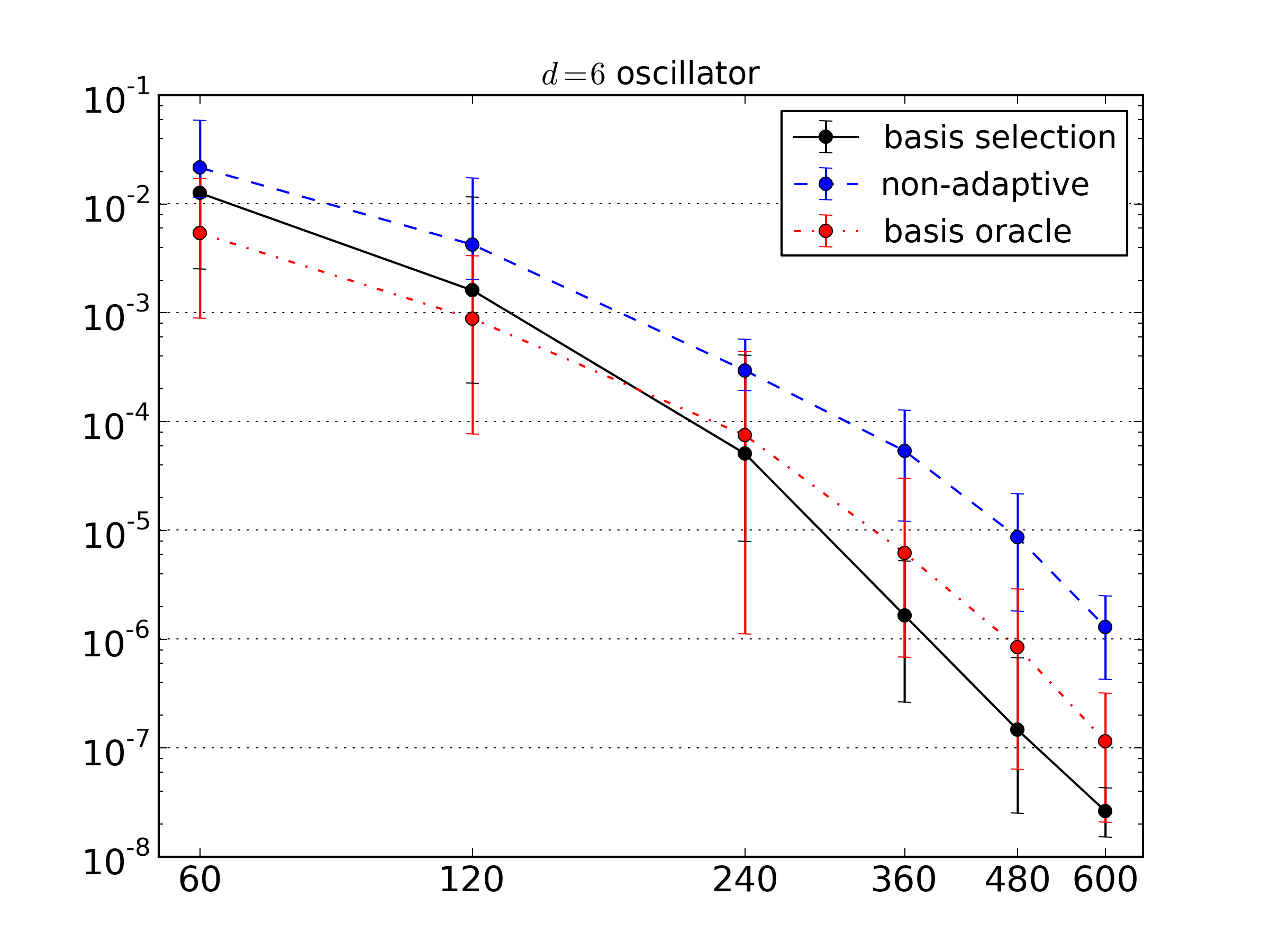}}\hspace{-0.95cm}
\subfloat[]{\includegraphics[width=0.7\textwidth]{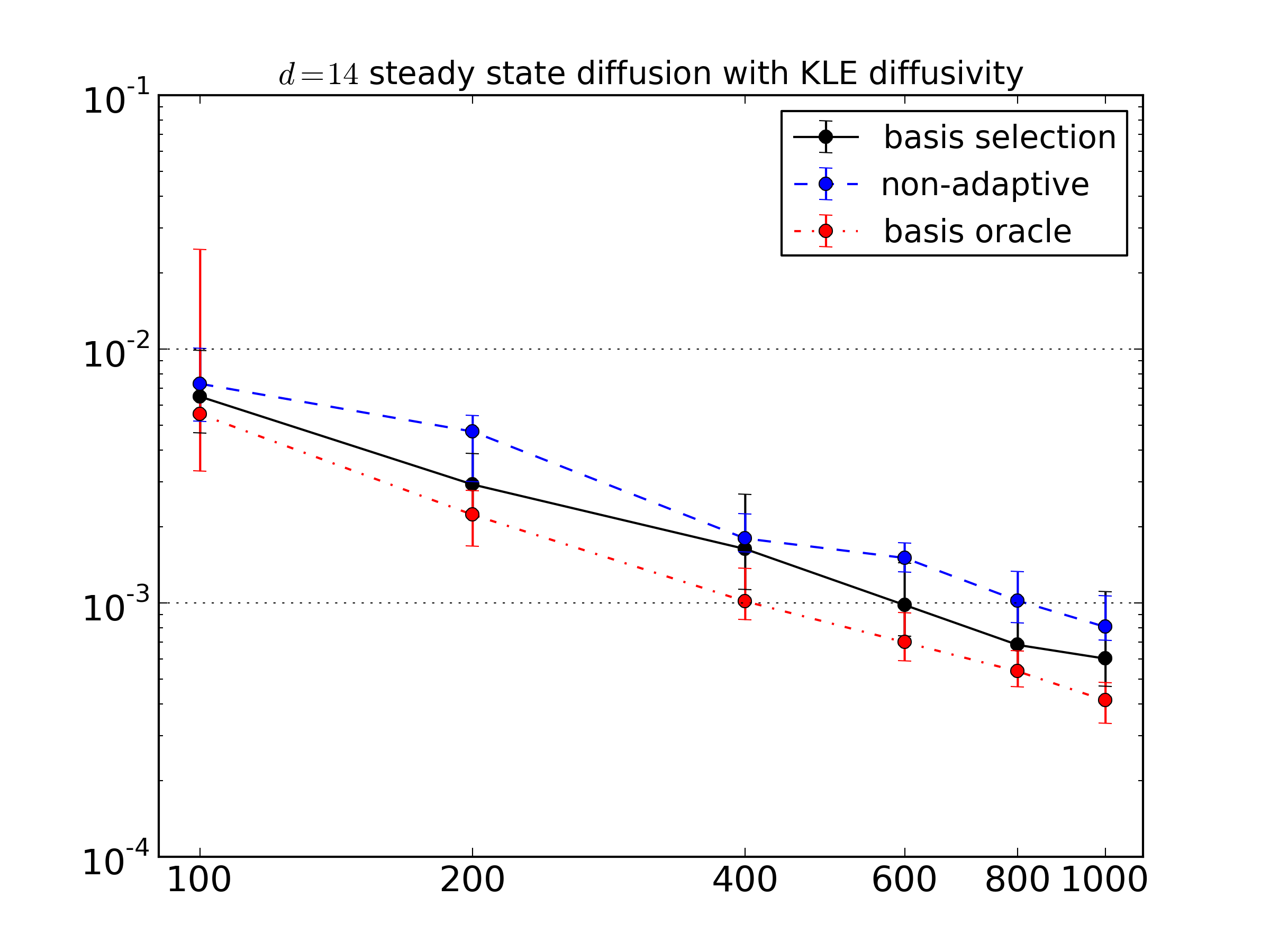}}
}
\makebox[1\linewidth][c]{%
\subfloat[]{
\includegraphics[width=0.7\textwidth]{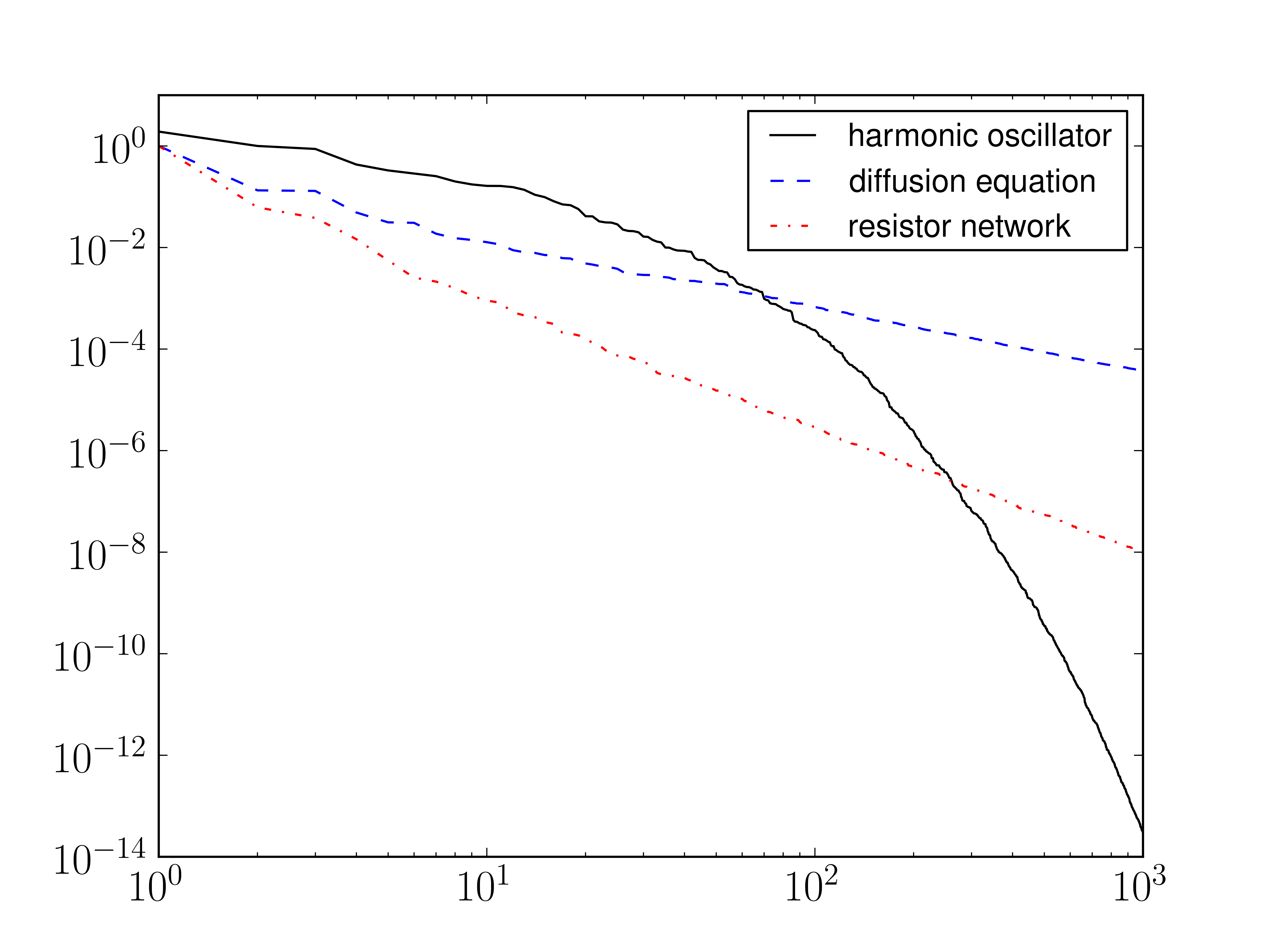}}
}
\caption{Convergence of the RMSE, with respect to increasing design size $M$, in the Legendre PCE approximation of (a) the harmonic oscillator~\eqref{eq:random-oscillator}
and (b) the solution $u(1/3,\brv)$ of the diffusion equation~\eqref{eq:hetrogeneous-diffusion}. 
The error bars represent the minimum and maximum error over the $20$ trials for each design size $M$.
(c) PCE coefficients, sorted by magnitude, for the harmonic oscillator, the diffusion equation and the resistor network.}
\label{fig:mixed-problems-coeff-decay}
% \label{fig:oscillator-d-6}
\end{figure}

Figure~\ref{fig:mixed-problems-coeff-decay} (a) depicts the error in the Legendre PCE for 
increasing design sizes $M$. The basis selection method clearly outperforms the non-adaptive approach and produces 
comparable results to the oracle. Again the improvement in performance is associated with a rapid decay of the exact PCE coefficients (see Figure~\ref{fig:mixed-problems-coeff-decay}).
% Again the increase in performance is associated with a high degree of anisotropy. Unlike the previous example 
% the importance of each variable on the model response does not decay with the dimension index. To provide insight
% into the nature of the anisotropy the total effect sensitivity indices~\eqref{eq:total-effect-indices} random oscillator are plotted in 
% Figure~\ref{fig:oscillator-d-6-total-effect-indices}. The majority of the variance in the function is attributed to
% the frequency $\omega$. The total effect indices can be extracted analytically from the PCE.
% \begin{figure}
% \centering
% \includegraphics[width=0.7\textwidth]{oscillator.png}
% \caption{Convergence of the RMSE, with respect to increasing design size $M$, in the Legendre PCE approximation of the harmonic oscillator~\eqref{eq:random-oscillator}. 
% The error bars represent the minimum and maximum error over the $20$ trials for each design size $M$.}
% \label{fig:oscillator-d-6}
% \end{figure}
% 
% \begin{figure}
% \centering
% \includegraphics[width=0.7\textwidth]{mixed-problems-coeff-decay}
% % \caption{PCE coefficients, sorted by magnitude, for the harmonic oscillator~\eqref{eq:random-oscillator}, the diffusion equation~\eqref{eq:hetrogeneous-diffusion} and the resistor network.}
% \caption{PCE coefficients, sorted by magnitude, for the harmonic oscillator, the diffusion equation and the resistor network.}
% \label{fig:mixed-problems-coeff-decay}
% \end{figure}

\subsection{Diffusion equation}
In this section, we consider the heterogeneous diffusion equation in one-spatial 
dimension subject to uncertainty in the diffusivity coefficient. This problem has been used as a benchmark
in other works~\cite{doostan2011,YangKarniadakis2013}.
Attention is restricted to one-dimensional physical space to avoid unnecessary 
complexity. The procedure described here can easily be extended to higher physical dimensions.
Consider the following problem with $d \ge 1$ random dimensions:
\begin{equation}\label{eq:hetrogeneous-diffusion}
-\frac{d}{dx}\left[a(x,\brv)\frac{du}{dx}(x,\brv)\right] = 1,\quad 
(x,\brv)\in(0,1)\times I_\brv
\end{equation}
subject to the physical boundary conditions
\begin{equation}
u(0,\brv)=0,\quad u(1,\brv)=0.
\end{equation}
Furthermore, assume that the random diffusivity satisfies
\begin{equation}\label{eq:diffusivityZ}
a(x,\brv)=\bar{a}+\sigma_a\sum_{k=1}^d\sqrt{\lambda_k}\phi_k(x)\rv_k,
\end{equation} 
where $\{\lambda_k\}_{k=1}^d$ and $\{\phi_k(x)\}_{k=1}^d$ are, respectively, 
the eigenvalues and eigenfunctions of the squared exponential covariance kernel 
\[
 C_a(x_1,x_2) = \exp\left[-\frac{(x_1-x_2)^2}{l_c^2}\right].
\]
The variability of the diffusivity field~\eqref{eq:diffusivityZ} is 
controlled by $\sigma_a$ and the correlation length
$l_c$ which determines the decay of the eigenvalues $\lambda_k$. Here we 
approximate the solution $u(1/3,\brv)$ with $\bar{a}=0.1$,
$d=14$, $\sigma_a=0.03$, $l_c=1/5$, while the uncertain inputs $\rv_k\in[-1,1]$, $k=1,\ldots,d$ are independent and uniformly distributed random variables.
% We also consider $d=40$, $\sigma_a=0.021$ and $l_c=1/14$. 
We solve the model~\eqref{eq:hetrogeneous-diffusion} using
quadratic finite elements with a high enough spatial resolution to neglect
discretization errors in our analysis. 

Figure~\ref{fig:mixed-problems-coeff-decay} (b) plots the error in Legendre PCE approximations built using
increasing design sizes $M$. There is negligible difference between
basis selection and the non-adaptive strategy, but there is also negligible difference between these methods and the oracle, indicating
there is not much improvement that can in principle be gained from basis selection. 
The negligible improvement is due to the fact that the `exact' PCE is not very compressible, as can be seen from Figure~\ref{fig:mixed-problems-coeff-decay}. The lack of compressibility means that many coefficients
are of similar magnitude and thus $\ell_1$-minimization in any form is not very effective.
%For the 40-dimensional example the dominant PCE coefficients correspond 
%to basis terms of degree three of less. Indeed the cross validation employed by the non-adaptive method chooses
%a degree three PCE as the final approximation. A similar argument can be made for the 14-dimensional case.
% \begin{figure}
% \centering
% \includegraphics[width=0.7\textwidth]{steady-heat-kle-d-14.png}
% %\includegraphics[width=0.45\textwidth]{steady-heat-kle-d-40.png}
% \caption{Convergence of the RMSE, with respect to increasing design size $M$, in the Legendre PCE approximation of the solution $u(1/3,\brv)$ of the diffusion equation~\eqref{eq:hetrogeneous-diffusion}.
% The error bars represent the minimum and maximum error over the $20$ trials for each design size $M$.}
% \label{fig:heat-eq-convergence-d-14}
% \end{figure}
% Note that when the 40 dimensional example was solved in~\cite{doostan2011}
% $\bset$ is chosen to contain all the terms of a 3rd-order Legendre PC expansion, i.e. $p = 3$ and the first 320 basis function from the $4$th-order chaos.
% Aside from the accuracy improvements show in the previous examples and the following, basis selection has the additional attraction of not having to resort
% to such ad-hoc techniques of choosing the basis $\bset$.

\subsection{Resistor network}
\begin{figure}
\includegraphics[width=0.95\textwidth]{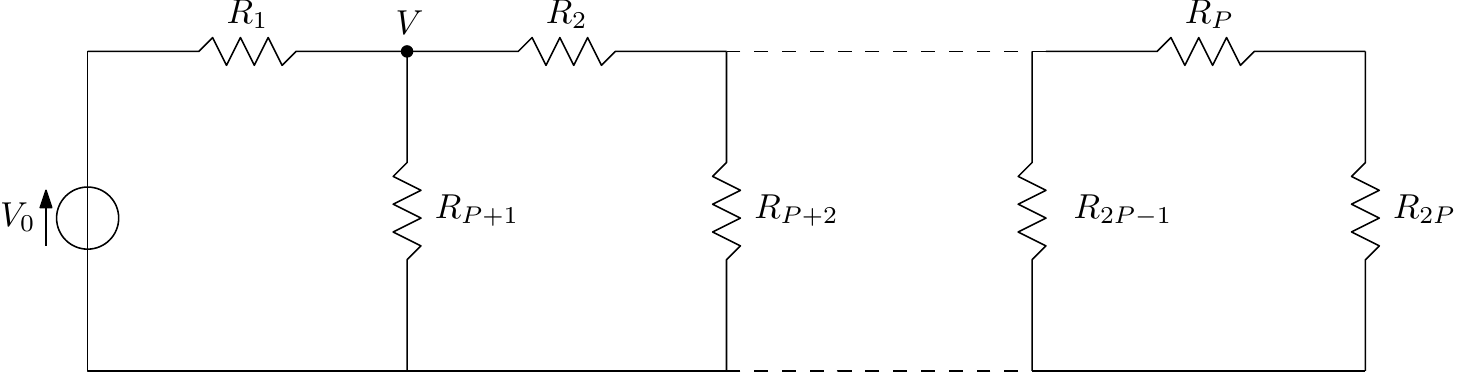}
\caption{Resistor network comprised of $d=2P$ resistances $R_i$, $i=1,\ldots,d$ of uncertain
ohmage and the network is driven by a voltage source providing a
known potential $V_0$. We are interested in determining voltage at $V$.}
\label{fig:resistor-network}
\end{figure}
As our last example, consider the electrical resistor network shown in Figure~\ref{fig:resistor-network}. 
The network is comprised of $d=2P$ resistances of uncertain
ohmage and the network is driven by a voltage source providing a
known potential $V_0$.  We are interested in determining how the voltage $V$ shown in 
the figure depends on the $d=2P$
resistances, which we take as random parameters uniformly distributed in the interval 
$\rv_k\in[1-\varepsilon, 1+\varepsilon]$, $k=1,\ldots,d$. This function is anisotropic. The effect of the resistors on the voltage
will decay with distance (in terms of the number of preceding resistors) from the point $V$.
In this example we set $d=40$ ($P=20$) and $d=80$ ($P=40$), take the maximum perturbation to be $\varepsilon = 0.1$ 
and set the reference potential $V_0 = 1$.

\begin{figure}
\centering
\makebox[1\linewidth][c]{%
\subfloat[]{\includegraphics[width=0.7\textwidth]{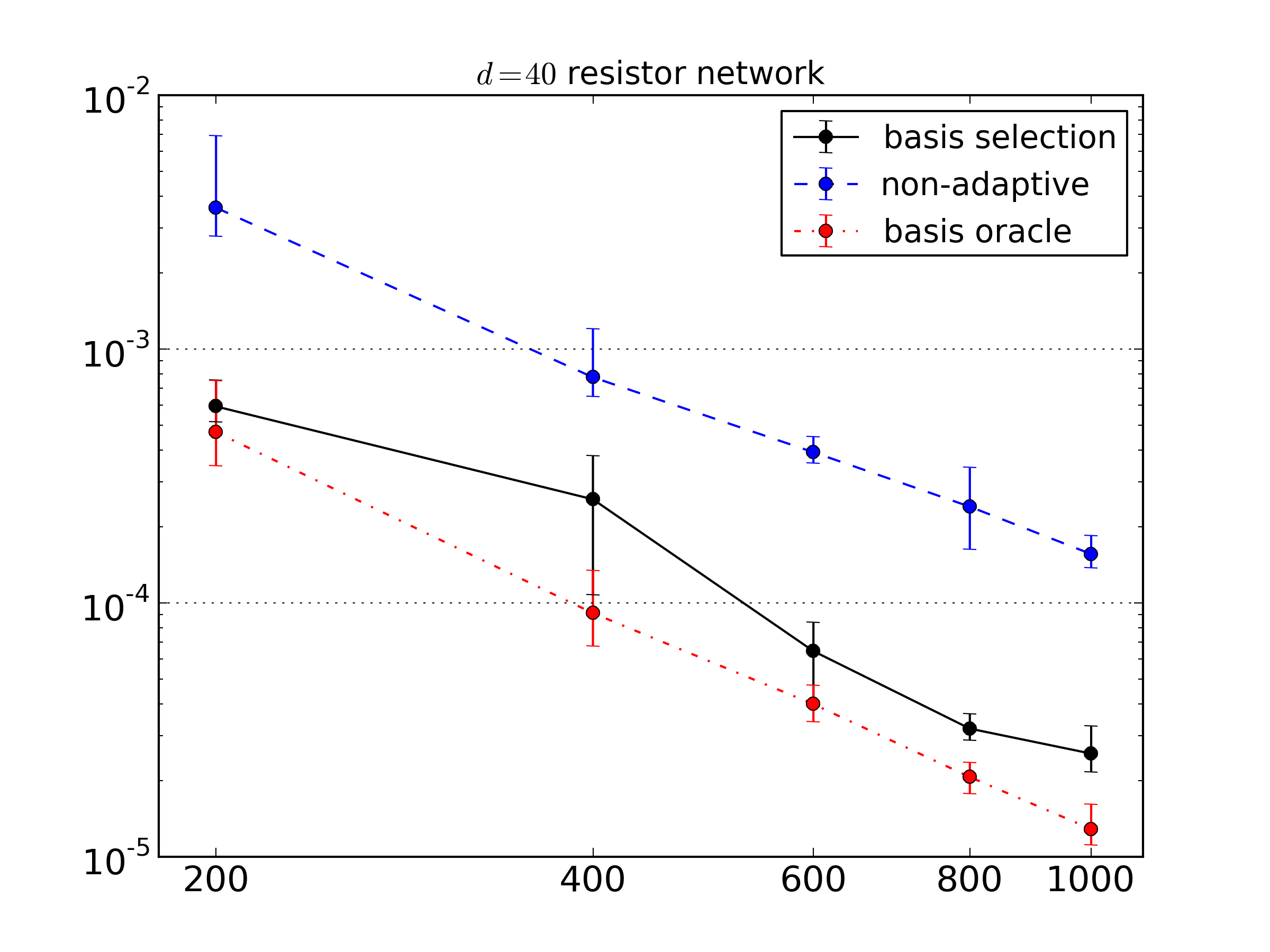}}
\hspace{-.95cm}
\subfloat[]{\includegraphics[width=0.7\textwidth]{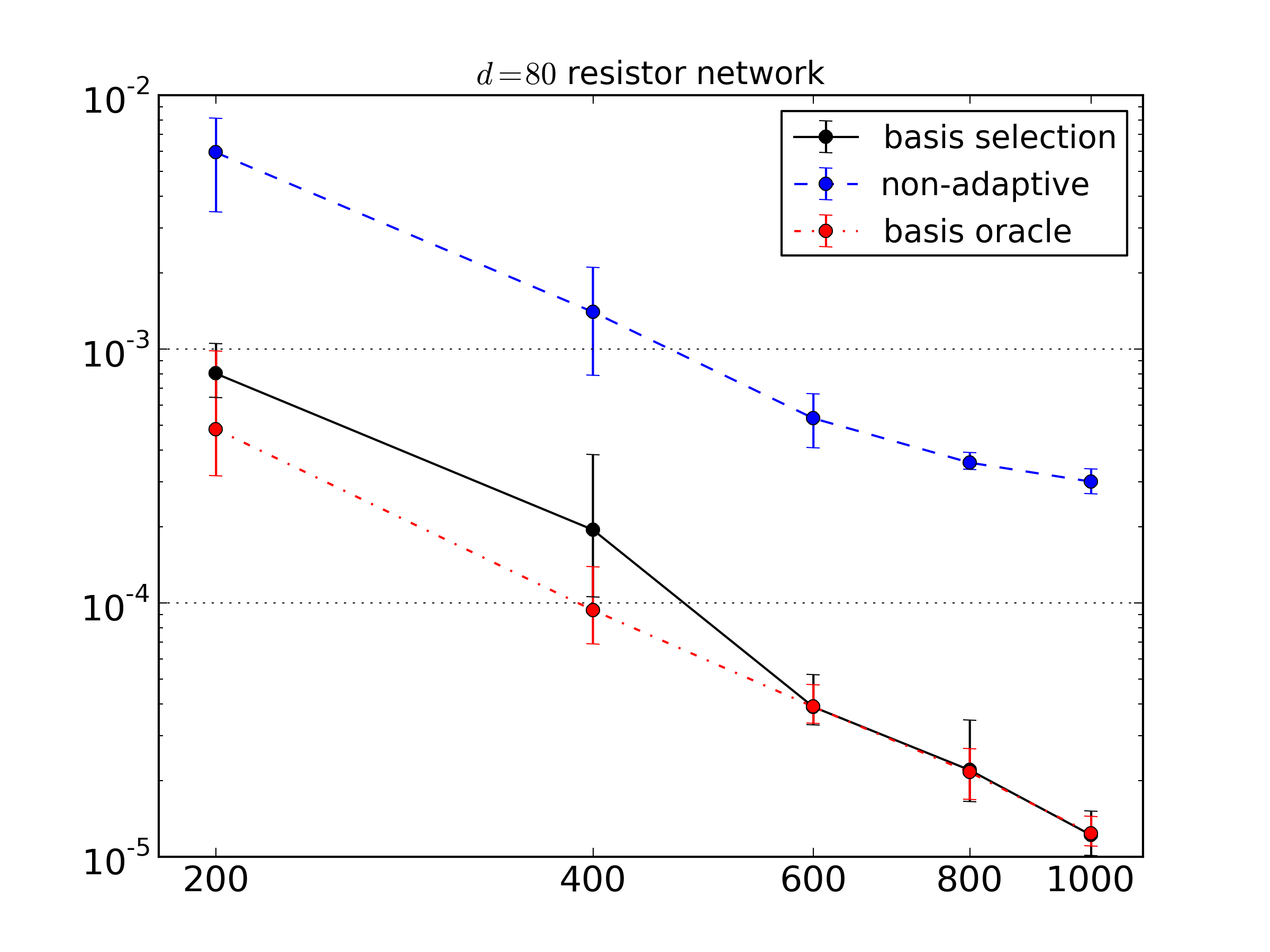}}
}
\caption{Convergence of the RMSE, with respect to increasing design size $M$, in the Legendre PCE approximation of the resistor network.
The error bars represent the minimum and maximum error over the $20$ trials for each design size $M$.}
\label{fig:resistor-d-40}
\end{figure}

Figure~\ref{fig:resistor-d-40}
shows the error in the Legendre PCE for increasing design sizes $M$. In both cases the basis selection method produces a PCE that is significantly 
more accurate than the PCE produced by the non-adaptive strategy. The basis selection method provides comparable results to the approximately optimal oracle.
% however, there is
% a degradation in accuracy with the smaller sample sizes, due to the inaccuracy of cross validation for small sample sizes. When the sample size is small
% splitting the data into 10 folds can significantly reduce the accuracy of the solution returned by OMP and, in turn, the ability to correctly
% choose the hyper-parameters $t$ and $\stoptol$ . As the sample size increases the effect of removing 10\% of the data becomes less significant.

\subsection{Gradient-enhanced $\ell_1$-minimization}
Typical $\ell_1$-minimization, when used for PCE approximation, attempts to find solutions to
\[
 \vand\coef \approx \boldf
\]
where denotes the Vandermonde matrix with entries
$\bPhi_{ij} = \phi_j(\brv_i),\quad i=1,\ldots,M,\; j=1,\ldots,N$. If gradients of the model $f$ with respect to the random variables $\brv$ are known, then
one can enhance the accuracy of the PCE by finding a solution to
\begin{equation*}
\vand_\partial\coef \approx \boldf_\partial
\quad\text{where}\quad
\vand_\partial=
\begin{bmatrix}
\vand \\ {\frac{\partial\vand}{\partial \rv_1}} \\ \vdots \\ {\frac{\partial\vand}{\partial \rv_d}}
\end{bmatrix},\quad
\boldf_\partial
\begin{bmatrix}
\boldf \\ {\frac{\partial\boldf}{\partial \rv_1}} \\ \vdots \\ {\frac{\partial\boldf}{\partial \rv_d}}
\end{bmatrix} 
\end{equation*}
and $({\frac{\partial\vand}{\partial \rv_n}})_{ij} = \frac{\partial \phi_j}{\rv_n}(\brv_i)$ and $(\frac{\partial\boldf}{\partial \rv_n})_i= \frac{\partial f}{\rv_n}(\brv_i)$,
$i=1,\ldots, M$, $j=1,\ldots,N$, $n=1,\ldots,d$. To find a solution we again use basis pursuit denoising and solve
\begin{equation*}
\label{eq:bpdn-grad}
\coef = \argmin_{\coef}\; \|\coef\|_1\quad \text{such that}\quad \|\bPhi_\partial\coef - \boldf_\partial\|_2 \le \stoptol
\end{equation*}
This gradient based formulation consists of $M(d+1)$ equations that match both function values and gradients, in comparison to~\eqref{eq:bpdn} which consists of only $M$ equations that match function values.

\begin{figure}
\centering
\makebox[1\linewidth][c]{%
\subfloat[]{\includegraphics[width=0.7\textwidth]{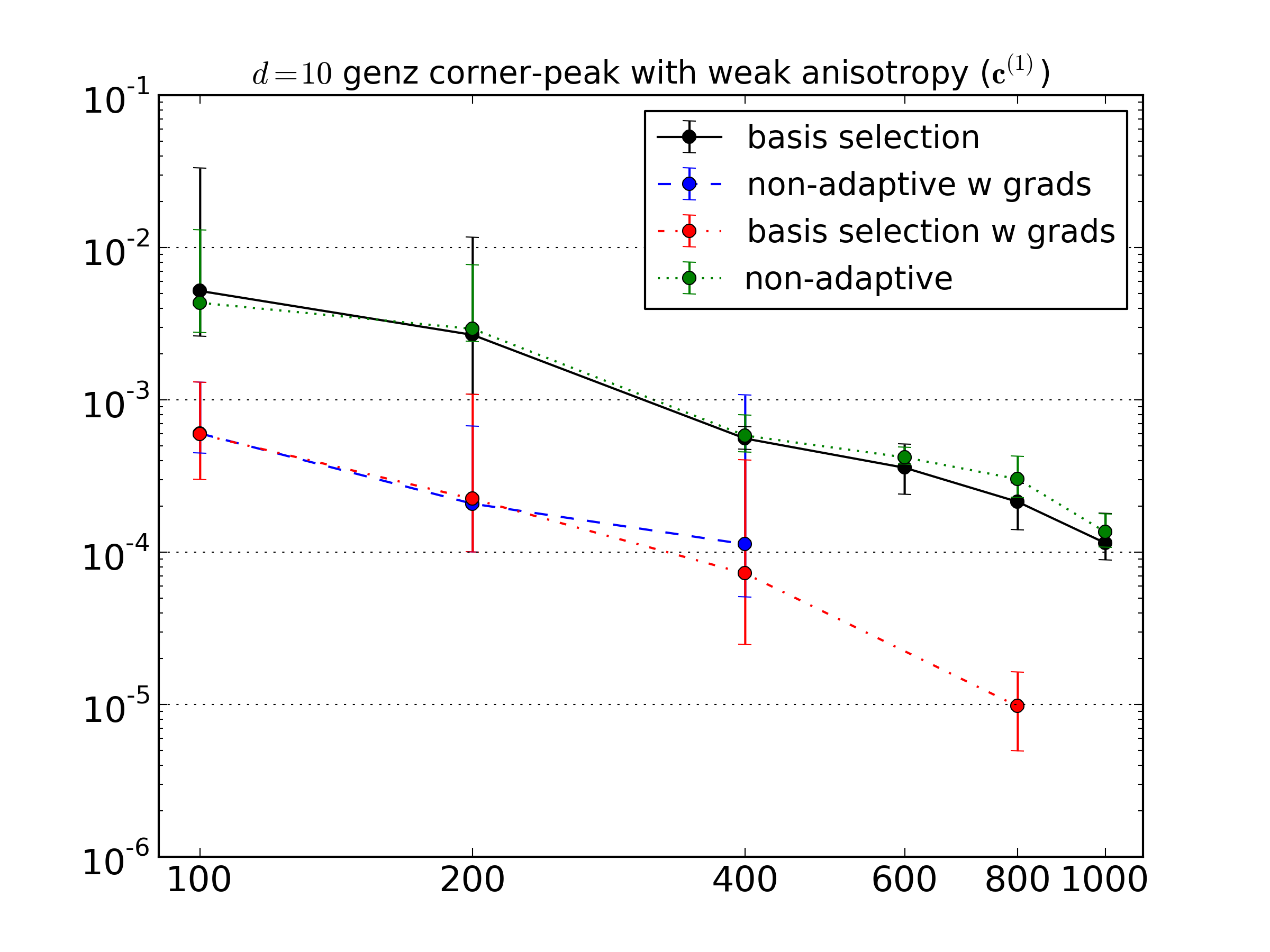}}
\hspace{-0.95cm}
\subfloat[]{\includegraphics[width=0.7\textwidth]{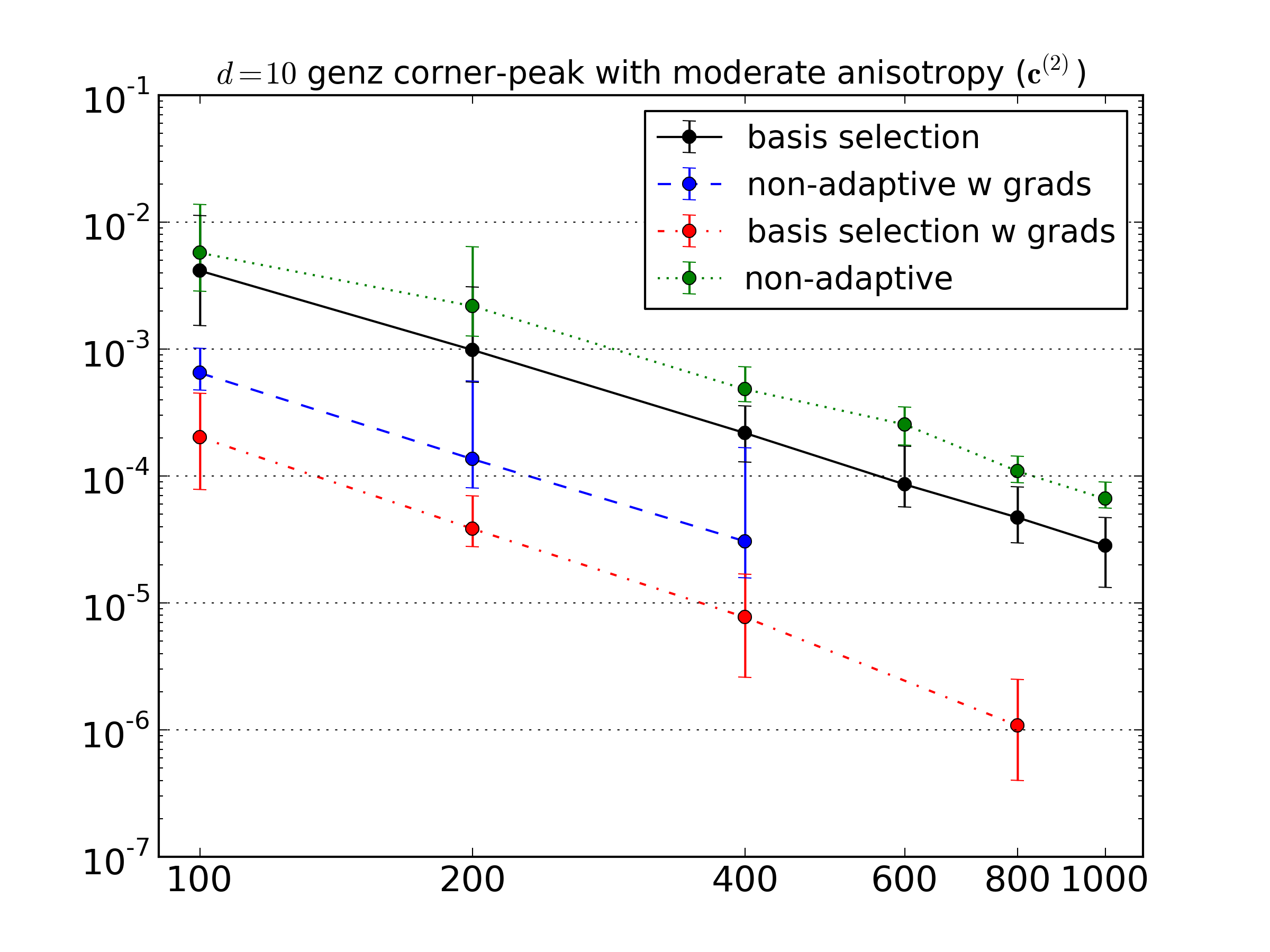}}
}
\makebox[1\linewidth][c]{
\subfloat[]{\includegraphics[width=0.7\textwidth]{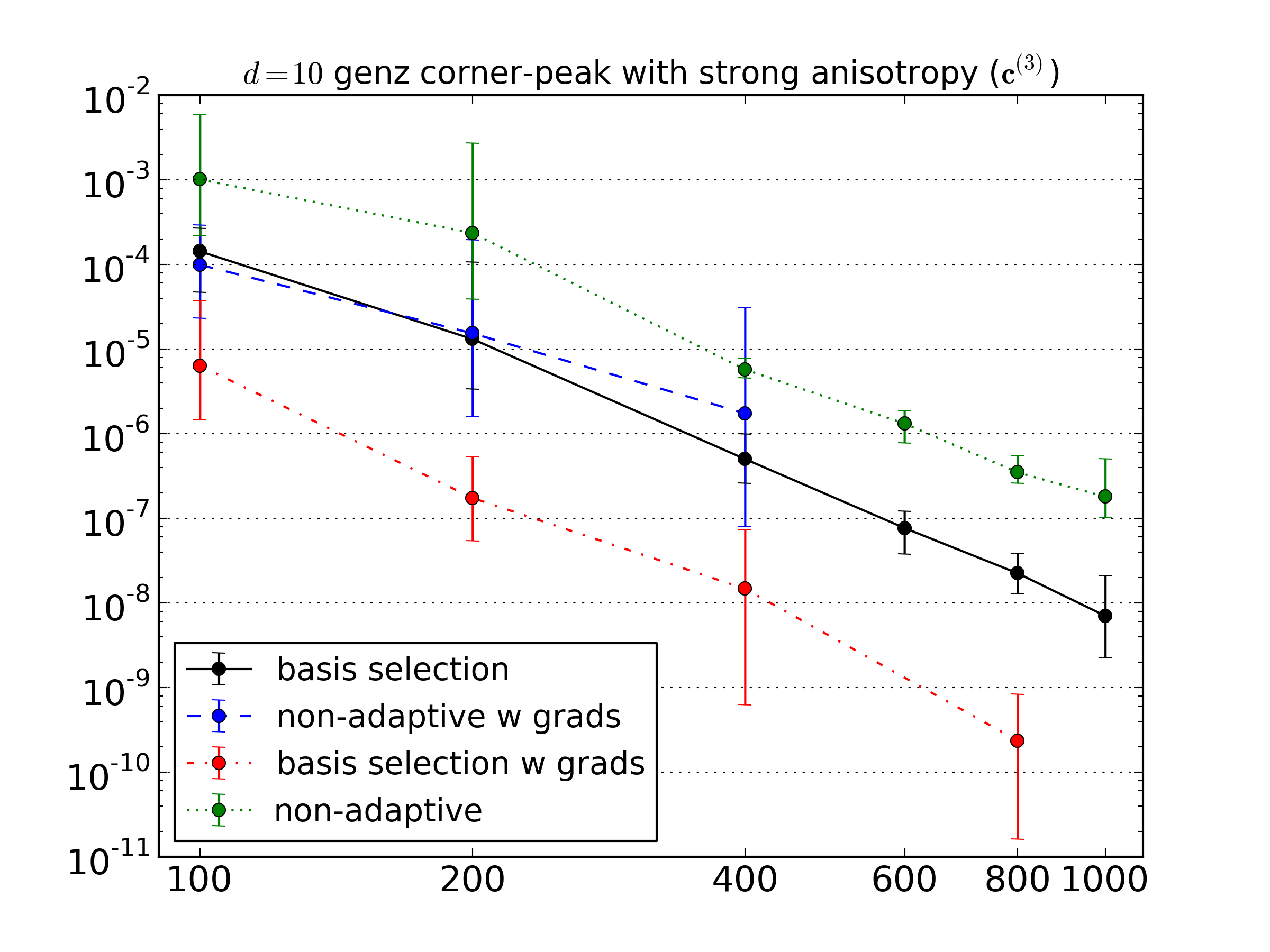}}
}
\caption{Convergence of the RMSE, with respect to increasing computation units, in the gradient-enhanced Legendre PCE approximation of the corner-peak function~\eqref{eq:genz-corner-peak}
for the three coefficient regimes: (a) $\bc^{(1)}$, (b) $\bc^{(2)}$, (c) $\bc^{(3)}$.
The error bars represent the minimum and maximum error over the $20$ trials for each design size $M$. We have assumed computing functional values at a design point is one unit and
computing all gradient components is another unit. For example, the gradient enhanced basis selection approximation (dot-dash line) with 50 design points costs 100 units (horizontal axis value).}
\label{fig:gradient-enhanced-genz}
\end{figure}

Figure~\ref{fig:gradient-enhanced-genz} demonstrates the utility of using gradient data to build PCE approximations of the corner-peak function~\eqref{eq:genz-corner-peak} with $d=10$. 
Unlike the previous figures in this paper, the horizontal axis is no longer the number of model
runs but rather the computational cost. We assume that running the model to only obtain function values costs one computational unit and running the model to obtain 
both function values and all gradients components requires two units. For example, adjoint methods for differential equations can be used to obtain all gradients at a
cost less than or equal to the cost of one forward model run.

Despite the extra computational cost required to obtain gradients, the use of gradients improves both the PCE resulting from both the non-adaptive and basis
selection methods. Similar to the results presented in Section~\ref{sec:genz-function-based-cs} the results shown here demonstrate that basis selection is more 
accurate than the non-adaptive strategy. Again the relative benefit is dependent on the rate of decay of the PCE coefficients. 

Basis selection is able to make effective use of gradient information. For a design $\bXi$ with $M$ samples, the size of the gradient enhanced Vandermonde matrix is $M(d+1)\times N$. 
We see that for a given accuracy gradient-based PCE requires a factor of $4$ fewer samples than the PCE based on the function values only. 
This is close to the optimal reduction factor of $d/2=5$ that can be obtained using gradients, 
assuming that each gradient component is as informative as a function value and the cost of computing function values with gradients is twice the 
cost of just computing function values.

% We remark that solving the gradient enhanced linear system is more computationally demanding for the same number of samples $M$ because the size of the linear system is $d+1$ times larger. 
% Thus the non-adaptive when used to compute PCE of approximately $100,000$ terms is much slower than the basis selection method, thus the convergence curves are not extended to the larger $M$
% used for the non-gradient enhanced approximations.

\subsection{Non-uniform model inputs}
\begin{figure}
\centering
\includegraphics[width=0.7\textwidth]{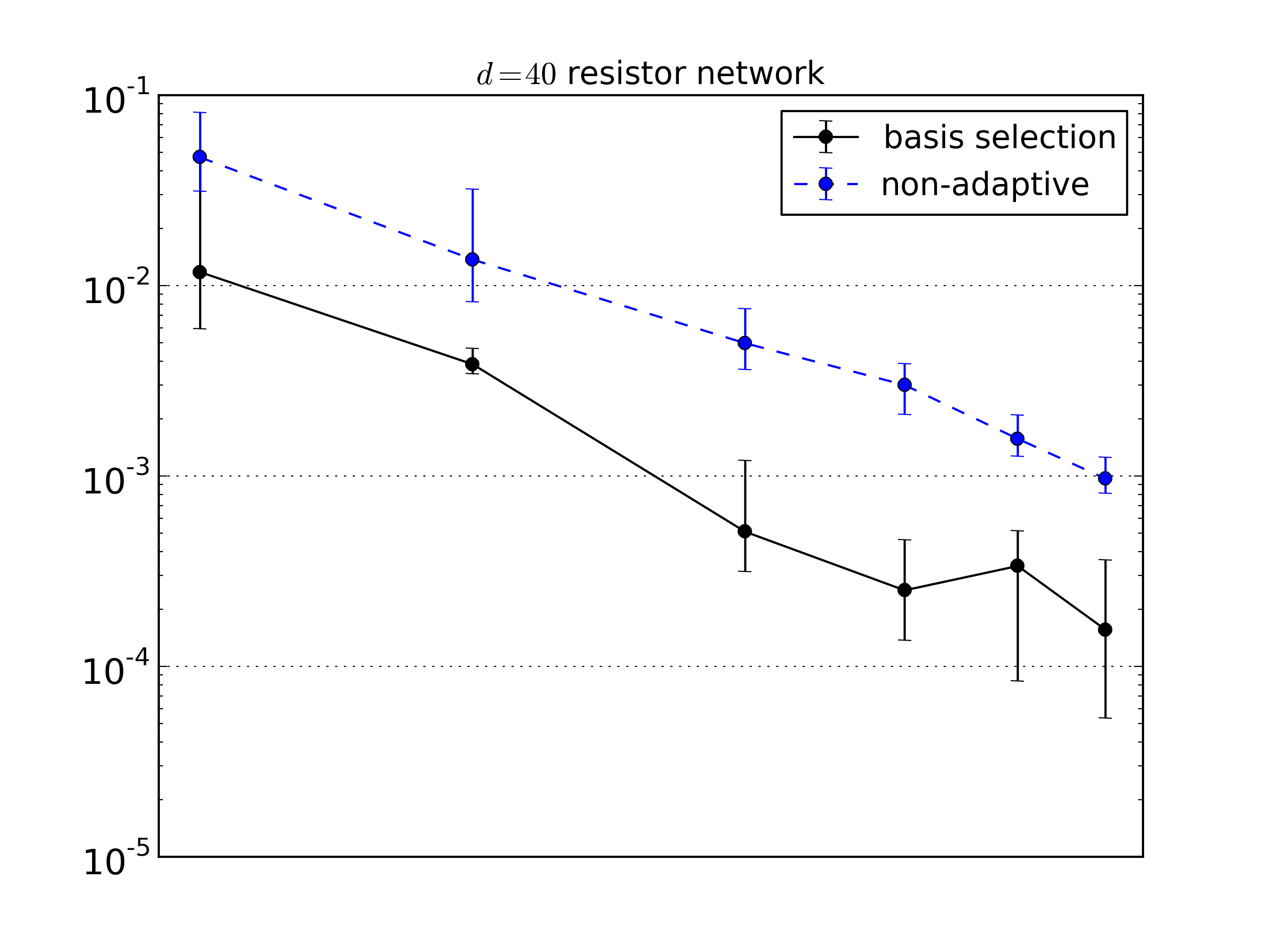}
\caption{Convergence of the RMSE, with respect to increasing design size $M$, in the Hermite PCE approximation of the $d=40$ resistor network. 
The error bars represent the minimum and maximum error over the $20$ trials for each design size $M$.}
\label{fig:normal-resitor-network-convergence}
\end{figure}

Throughout this paper, we have discussed basis selection when applied to Legendre polynomials and uniform variables. 
However, basis selection can also be applied to other variable/polynomial 
combinations. Let us once again consider the resistor network, but now let $\brv$ be Gaussian variables with mean $1.0$ and standard 
deviation $0.005.$\footnote{The standard 
deviation is made sufficiently small to make the chance of negative resistances practically zero.} We now draw random samples
from the aforementioned Gaussian distribution to form $\bXi$ and run the model at each
sample to obtain $f(\bXi)$. 
Figure~\ref{fig:normal-resitor-network-convergence} demonstrates that the advantages of basis selection are also present when 
we compute PCE approximations with non-uniform random variables. 

Note that linear systems based upon Hermite polynomials suffer from poor numerical conditioning as the number of samples $M$ increases. 
This poor conditioning causes the non-monotone convergence shown.
Development of sampling and pre-conditioning strategies for normal variables is an important area of future research, but is beyond the scope of this paper.

\section{Conclusions}
\label{sec:conclusions}
In this paper we present a basis selection method that can be used with
$\ell_1$-minimization to adaptively determine the large coefficients of 
polynomial chaos expansions (PCE). The method attempts to identify structure in
the coefficients of a PCE and only applies $\ell_1$-minimization to those terms
believed to have large coefficients. The adaptive construction produces anisotropic basis sets that have more terms in important dimensions and limits the number of unimportant terms
which increase mutual coherence and thus degrade the performance of $\ell_1$-minimization. The basis selection
method produces, for a given computational budget, a more 
accurate PCE than would be obtained if the basis is fixed a priori.
The important features and the accuracy of basis selection are demonstrated with a number of numerical examples. 
Specifically we show that in high dimensions, for which high-order
total-degree PCE bases are infeasible, basis selection allows accurate
identification of high-order terms that cannot be captured by low-order total-degree expansions. 
We demonstrate that even for lower dimensional problems, for
which high-order total-degree PCE bases are feasible, 
basis selection still produces more accurate results than basis
sets that are fixed priori. Finally, we demonstrate that basis selection can effectively leverage function gradients
and be applied to PCE of non-uniform random variables.

\bibliographystyle{plain}
%\bibliography{../../bibliography/gaussian-processes,../../bibliography/EpistemicUncertainty,../../bibliography/Collocation,../../bibliography/SparseGrid,../../bibliography/Misc,../../bibliography/GPC,../../bibliography/JDJakeman,../../bibliography/DiscontinuityDetection,../../bibliography/compressed-sensing,../../bibliography/cross-validation,../../bibliography/Cubature}
\bibliography{references}
\end{document}